\documentclass[a4paper,12pt,reqno]{amsart}

\usepackage[letterpaper, hmargin=1in, top=1in, bottom=1.2in, footskip=0.7in]{geometry}
\usepackage{rotating,pdflscape,xcolor,amssymb,subfigure,psfrag,amsmath,eufrak,bbm,epsfig,amsthm,mathtools,fancyhdr,graphicx}
\definecolor{vdarkred}{rgb}{0.6,0,0.2}
\definecolor{vdarkblue}{rgb}{0,0.2,0.6}
\usepackage[pdftex, colorlinks, linkcolor=vdarkblue,citecolor=vdarkred]{hyperref}
\usepackage{a4wide,amscd}
\usepackage{changepage}
\usepackage{centernot}

\usepackage{tikz} 
\usepackage{pgfplots}

\usepackage{listings}
\usepackage[small]{caption}

\newcommand{\ii}{\mathrm{i}}
\newcommand{\me}{\mathrm{e}}

\newcommand{\cA}{\mathcal{A}}\newcommand{\cB}{\mathcal{B}}
\newcommand{\cC}{\mathcal{C}}

\newcommand{\cN}{\mathcal{N}}

\newcommand{\bC}{\mathbf{C}}
\newcommand{\bE}{{\mathbf{C}_{XY}}}
\newcommand{\bERF}{{\mathbf{C}_{YX}}}

\newcommand{\bERFoos}{\mathbf{C}_{YX}^{\op{out\,of\,sample}}}
\newcommand{\bF}{\mathbf{F}}
\newcommand{\bG}{\mathbf{G}}

\newcommand{\bU}{\mathbf{U}}\newcommand{\bV}{\mathbf{V}}
\newcommand{\bX}{\mathbf{X}}
\newcommand{\bY}{\mathbf{Y}}\newcommand{\bZ}{\mathbf{Z}}

\newcommand{\bee}{\mathbf{e}}

\newcommand{\bu}{\mathbf{u}}\newcommand{\bv}{\mathbf{v}}

\newcommand{\wta}{\widetilde{a}}\newcommand{\wtb}{\widetilde{b}}

\newcommand{\wtA}{\widetilde{A}}\newcommand{\wtB}{\widetilde{B}}

\newcommand{\wtG}{\widetilde{\bG}}

\newcommand{\ka}{\kappa} 

\newcommand{\si}{\sigma}
\newcommand{\lam}{\lambda}

\newcommand{\vfi}{\varphi}
\newcommand{\al}{\alpha}
\newcommand{\tta}{\theta}
\newcommand{\Tta}{\Theta}
\newcommand{\Si}{\Sigma}
\newcommand{\bet}{\beta}
\newcommand{\del}{\delta}

\newcommand{\ld}{\ldots}
\newcommand{\cdote}{\,\cdot\,}

\newcommand{\beg}{\begin}
\newcommand{\en}{\end}

\renewcommand{\Im}{\mathfrak{Im}}

\newcommand{\trm}{\textrm}

\newcommand{\bgt}{\begin{itemize}}
\newcommand{\ent}{\end{itemize}}
\newcommand{\ite}{\item}
\newcommand{\eqre}{\eqref}
\newcommand{\re}{\ref}
\newcommand{\la}{\label}

\newcommand{\brem}{\begin{rmk}}
\newcommand{\erem}{\end{rmk}}
\newcommand{\blem}{\begin{lem}}
\newcommand{\elem}{\end{lem}}
\newcommand{\bcor}{\begin{cor}}
\newcommand{\ecor}{\end{cor}}
\newcommand{\bTh}{\begin{Th}}
\newcommand{\eTh}{\end{Th}}
\newcommand{\bpropo}{\begin{propo}}
\newcommand{\epropo}{\end{propo}}

\newcommand{\op}{\operatorname}

\newcommand{\Cov}{\operatorname{Cov}}

\newcommand{\diag}{\operatorname{diag}}

\newcommand{\Tr}{\operatorname{Tr}}

\newcommand{\ud}{\mathrm{d}}

\newcommand{\mE}{{m}_{\bE,\cC}}

\newcommand{\E}{\op{\mathbb{E}}}

\newcommand{\R}{\mathbb{R}}
\newcommand{\C}{\mathbb{C}}

\newcommand{\p}{\mathbb{P}}

\newcommand{\f}{\frac}
\newcommand{\ff}{\frac{1}}
\newcommand{\lf}{\left}
\newcommand{\ri}{\right}

\newcommand{\st}{such that }

\newcommand{\ti}{\times}

\newcommand{\mc}{\mathcal}

\newcommand{\eps}{\varepsilon}

\newcommand{\bck}{\backslash}

\newcommand{\bbm}{\begin{bmatrix}}
\newcommand{\ebm}{\end{bmatrix}}
\newcommand{\bes}{\begin{equation*}}
\newcommand{\ees}{\end{equation*}}
\newcommand{\be}{\begin{equation}}
\newcommand{\ee}{\end{equation}}
\newcommand{\beqy}{\begin{eqnarray}}
\newcommand{\eeqy}{\end{eqnarray}}
\newcommand{\beq}{\begin{eqnarray*}}
\newcommand{\eeq}{\end{eqnarray*}}

\newcommand{\ie}{i.e. }

\newcommand{\bpm}{\begin{pmatrix}}
\newcommand{\epm}{\end{pmatrix}}

\newcommand{\cd}{\cdots}

\newcommand{\bpr}{\beg{proof}}
\newcommand{\epr}{\en{proof}}

\newcommand{\nober}{\nonumber}

\newcommand{\pa}{\partial}

\newcommand{\FrameTable}[2]{#1\scriptsize{ $\pm$#2}}

\newcommand{\whcC}{\widehat{\cC}}
\newcommand{\AND}{\qquad\trm{ and }\qquad}

\newcommand{\ucleaned}{^{\op{cleaned}}}
\newcommand{\uAlgo}{^{\op{Algo}}}
\newcommand{\uLedoitPeche}{^{\op{Ledoit-Peche}}}
\newcommand{\ucleanedalgo}{^{\op{cleaned,\,algo}}}
\newcommand{\utrue}{^{\op{true}}}
\newcommand{\Xtrue}{^{\op{true},X}}
\newcommand{\Ytrue}{^{\op{true},Y}}
\newcommand{\uclean}{^{\op{clean}}}
\newcommand{\oos}{_{\op{oos}}}

\newcommand{\meanoosovl}{\op{OVL}_{oos}}
\newcommand{\meanisovl}{\op{OVL}_{is}}

 \newcommand{\theo}[1]{Theorem \re{#1}}
\newcommand{\co}[1]{Corollary \re{#1}}
\newcommand{\lemm}[1]{Lemma \re{#1}}
\newcommand{\rmq}[1]{Remark \re{#1}}
\newcommand{\prop}[1]{Proposition \re{#1}}

\newcommand{\fig}[1]{Figure \re{#1}}

\allowdisplaybreaks
\newtheorem{Th}{Theorem}[section]

\newtheorem{propo}[Th]{Proposition}

\newtheorem{lem}[Th]{Lemma}

\newtheorem{cor}[Th]{Corollary}

\theoremstyle{definition}
\newtheorem{rmk}[Th]{Remark}

\long\def\symbolfootnote[#1]#2{\begingroup
\def\thefootnote{\fnsymbol{footnote}}\footnote[#1]{#2}\endgroup}

\definecolor{gris25}{gray}{0.75}
\newcounter{algo}
\newenvironment{algo}[1]
{
\refstepcounter{algo}
\small \sf \medskip
\noindent\colorbox{gris25}{
\makebox[\textwidth][c]{\bfseries #1}}\bigskip 
}
{\smallskip
\noindent\colorbox{gris25}{\makebox[\textwidth][c]{\hspace{1cm}}}
\rm \normalsize \smallskip}


\author{Florent Benaych-Georges}
\author{Jean-Philippe Bouchaud}
\author{Marc Potters}

\address[FBG, JPB, MP]{CFM, 23 rue de l'Universit\'e, 75007 Paris, France}
\email[FBG]{florent.benaych-georges@cfm.fr}
\email[JPB]{jean-philippe.bouchaud@cfm.fr}
\email[MP]{marc.potters@cfm.fr}

\date{\today}
\title{Optimal cleaning for singular values of cross-covariance matrices}
\begin{document}

\maketitle
 \begin{abstract}We give a new algorithm for the estimation of the cross-covariance matrix $\mathbb{E} XY'$ of two large dimensional signals $X\in\mathbb{R}^n$, $Y\in \mathbb{R}^p$ in the context where the number $T$ of observations of the pair $(X,Y)$ is large but $n/T$ and $p/T$ are  not supposed to be small.  In the asymptotic regime where $n,p,T$ are large, with high probability, this algorithm  is optimal for the Frobenius norm  among \emph{rotationally invariant estimators}, i.e. estimators derived from the \emph{empirical estimator} by \emph{cleaning} the singular values, while letting singular vectors unchanged. 
 \en{abstract}
    \tableofcontents

  \section{Introduction}

   \subsection{Context} 
    In high-dimensional statistics, it is well known that the classical \emph{empirical estimator}  (\ie the one based on  an average over the sample)   
	has little efficiency   when the sample size is not much larger than the dimension of the object we want to estimate. 
	For example, the spectrum of the  empirical covariance matrix of a sample of $T$ independent observations of an $n$-dimensional  Gaussian signal with covariance $I_n$ 
	is not concentrated in the neighborhood of $1$ when $T$  has the same order as $n$, but distributed according to the Marchenko-Pastur law with parameter $n/T$.    
	In the same way, for $(X(t),Y(t))_{t=1,\ldots, T}$ a sample of observations of  a pair  $(X,Y)\in \R^n\ti \R^p$ of random vectors,
	the singular values of the empirical estimator \be\la{EmpEstIntro}\qquad \qquad \bE:=\ff{T}\sum_{t=1}^T X(t)Y(t)' \qquad\text{with $Y(t)':=$ transpose of the column $Y(t)$}\ee of the true cross-covariance matrix 
	are not distributed as the singular values of the 
    true cross-covariance matrix  when $T$ is  not large with respect to $n$ and $p$ (see Figure \re{Fig_CFM_talk}, where we plot 
	both the true singular values density and the histogram of the empirical singular values). 
	\newcommand{\scaleCC}{0.5}
\begin{figure}[ht]
\centering
\includegraphics[scale=\scaleCC]{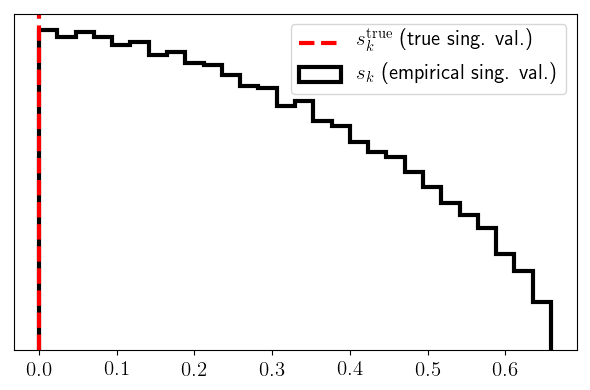}
\includegraphics[scale=\scaleCC]{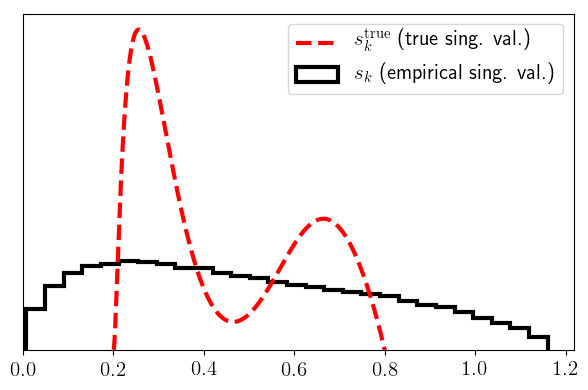}
\caption{Singular values of  $\ff{T}\sum_tX(t)Y(t)'$   vs true singular values. 
\textbf{Left:}   $\protect\bpm X\\ Y\protect\epm
 \sim\cN(0,I_{n+p})$.  \textbf{Right:} 
  $\protect\bpm X\\ Y\protect\epm\sim\cN(0,\Si)$ for $\Si=\protect\bpm I_n & \cC\\ \cC'&I_p\protect\epm$ with $\cC$ having singular values with density given by the red dashed curve. In both cases, $T/n=T/p=10$ and $T=25000$. The total lack of fit of the red curve by the histogram on the right and the spread between the true value $0$ and most of the 
 histogram on the left show that the empirical estimator works poorly (even though $T$ is 10 times larger than $n$ and $p$).}\label{Fig_CFM_talk}
\end{figure}

In the case of covariance estimation, a problem of interest in finance \cite{LCBP,LedoitWolf20042,BBPoRisk,BBPo},  several methods have been developed to circumvent these difficulties and improve the empirical estimator, 
based on regularization \cite{ElKaroui,bose,ElKaroui3},   shrinkage \cite{LedoitWolf2004,RIE0,LedoitWolf2012,BBPoRisk,BBPo}, 
specific sparsity or low-rank assumptions on the true covariance matrix \cite{ElKaroui2,klopp,KLT1,klopp2}, 
robust statistics \cite{Robust1,Robust2} or fixed-point analysis \cite{popescu}.

 However, the problem of the estimation of cross-covariance matrices has, to our knowledge, not been addressed so far, despite 
 its numerous applications in various fields (see e.g. \cite{CorrBouchaud}, where the null model is studied).
 
 Of course, cross-covariance estimation can formally be considered as a sub-problem of covariance estimation, as any pair  of  random vectors  $(X,Y)\in \R^n\ti \R^p$ 
 can be concatenated in a vector $Z=\bpm X\\ Y\epm\in \R^{n+p}$ whose covariance matrix has upper-right corner the cross-covariance of 
 $X$ and $Y$.  
 The problem with this idea is  that  the above covariance estimation methods rely on prior information on the structure of the covariance matrix: some of them are based on the hypothesis that the covariance matrix   of $Z$ is   sparse or low-rank or essentially supported by a neighborhood of its diagonal and some others   work in  the Bayesian framework  where  the true covariance matrix of $Z$
 has been chosen at random with a prior distribution that is invariant under the action of the orthogonal group by conjugation (\emph{rotationally invariant estimators} \cite{RIE0,LedoitWolf2012}), which implies that the entries of $Z$ can naturally be blended in linear combinations. This clearly does not make sense when  $X$ and $Y$ are of different nature, 
   for example if $X$ contains commodity price returns and $Y$, say, weather data\footnote{The prices of lots of commodities (e.g. energy, agricultural products) indeed exhibit  strong correlations with the weather.}. However, an
   analogue notion exists for cross-covariance matrices, that we also call 
 \emph{rotationally invariant estimators}: estimators based on the  empirical estimator from \eqref{EmpEstIntro}, modifying (we say \emph{cleaning}) its singular values, 
 but letting its singular vectors unchanged, 
 \ie estimators relevant to  the Bayesian framework  where  the  true 
 cross-covariance matrix has been chosen at random, with a prior distribution that is invariant under the actions of the orthogonal groups by multiplication on the left and on the right. 

 \subsection{Contents of the paper}\la{section_optimality}
 \subsubsection{Purpose}
   The purpose of this text is precisely to compute the optimal rotationally invariant estimator for the true cross-covariance   in the regime 
   where we have at disposal a large number $T$ of observations  of the pair $(X,Y)$, but where $n/T$ and $p/T$ are not supposed to be small. 
    It is \emph{optimal} in the sense that for Gaussian data, for $$\cC:=\E XY'$$ the true cross-covariance of $X$ and $Y$,  it is the solution of 
   \be\la{intro_optimal_eq}\op{argmin}_{\op{estimators}}
   \| \op{Estimator}-\cC\|_{\text{F}}\ee
   among the estimators   whose singular vectors 
  are those of  the \emph{empirical estimator} $\bE$ given at \eqre{EmpEstIntro} above. Here,   $\|\cdot\|_{\text{F}}$  denotes the \emph{Frobenius norm} ,  \ie   the standard Euclidean norm on matrices:   \be\la{opt_problem00}\|M\|_{\text{F}}:=\sqrt{\Tr MM'}.\ee

  Let us introduce the SVD of the empirical estimator $\bE$ from \eqre{EmpEstIntro}: $$\bE=\sum_ks_k\bu_k\bv_k',$$ with $s_k$ the singular values and $\bu_k$ (resp. $\bv_k$) the left (resp. right) singular vectors.  One easily gets  (see \eqre{0507181} below) that optimality rewrites 
  \be\la{def:sikskIntro} \op{Estimator}=\sum_ks_k\ucleaned\bu_k\bv_k'\quad\text{ with} \quad s_k\ucleaned=\bu_k'\cC\bv_k.\ee
  The numbers $$\bu_k'\cC\bv_k,$$  called  \emph{oracle estimates}, are of course unknown, and the main problem is to estimate them.  Rather than computing them directly, we shall introduce, at \eqre{Notation:Intro250700AA},   a function $L(z)$ ($z\in \C\bck\R$), which allows to estimate them (at least their weighted averages, which is enough for our purpose). The function $L$ is  called the \emph{oracle function} for a reason explained below (see \prop{propo:sol_optim} and right above it) and is   estimated in terms of observable variables   (Theorems \re{th:main} and \re{th:main2}).
  
  Note that the \emph{optimality of the matrix estimator} does not imply that as the dimensions tend to infinity, the numbers $s_k\ucleaned-s_k$ tend to zero  (see \eqre{rough_estimate_cleaned_true_emp2} and Proposition \re{prop_unbiased_estimator}). The reason is that  the optimal \emph{matrix estimator}   keeping \emph{empirical singular vectors unchanged}  takes into account the fact that  these vectors are noisy versions of the true ones, hence reduces their weights by shrinking the singular values.   
   
\subsubsection{Main contributions}
The main contributions of the paper are the following ones:
  \bgt \ite In Proposition \re{propo:sol_optim}, we express our \emph{oracle estimates} $\bu_k'\cC\bv_k$ in terms of the observable function $G(z)$ and of the unobservable function (called  \emph{oracle function} ) $L(z)$.\\
  \ite The main achievement of the paper is then to provide an approximation of our  \emph{oracle function} $L(z)$ defined at \eqre{Notation:Intro250700AA} in terms of observable variables   (Theorems \re{th:main} and \re{th:main2}).   The proofs, given in Section \re{sec:proofs},  are based on classical concentration results (Proposition \re{247153}) and on long  computations starting from the Stein formula for     Gaussian random vectors
 (Proposition \ref{SteinMultidim}). This approach is completely different from the one of the proofs of the Ledoit-P\'ech\'e paper \cite{RIE0} (where an analogue estimator for \emph{covariance} matrices was proposed) and, incidentally, allows to recover the main formula of \cite{RIE0} very directly (see \cite{FloRIE}). 
 \beg{rmk}One of the advantages of the proofs of Ledoit-P\'ech\'e's paper is that they do not rely on any Gaussian hypothesis. One  can then wonder whether the method of the present paper (and its main results) could be extended beyond  the Gaussian framework. It happens that the Stein formula can be generalized, with an error term,   beyond Gaussian variables (see e.g. \cite[Lem. 1.13.9]{FloAntti}) and that even though we did not include it here for brevity, we checked that the main lines of the proof still work under much more general assumptions (including of course the fact that $Z$ is centered, has covariance $\Si$ and moments up to a reasonable order). \en{rmk}
  \ite We provide precise, presumably close to sharp,  error terms for the approximation of the \emph{oracle function},  which, contrarily to \cite{RIE0}, do not rely on convergence hypotheses  for empirical spectral distributions and are controlled by very few quantities.\\
  \ite In   Section \re{sec:numsim}, devoted to numerical simulations, we assess  and illustrate algorithm accuracy by comparing it with the empirical estimator and the upper-right corner of  Ledoit-P\'ech\'e's estimator for a quite diversified set of models, which is not at all restricted to the invariance class this estimator was thought for (\ie the one described \eqre{1301181}). For all the  models we simulate, our algorithm outperforms (most times by far) both other algorithms in the sense that its output is closer to the true cross-covariance matrix, for the Frobenius norm as well as for the operator norm (see Tables \re{ref:my_table} and \re{ref:my_table_op}). The python code for the numerical simulations of this paper is available at \url{https://github.com/CFMTech/Optimal_cleaning_for_singular_values_of_cross-covariance_matrices}\\
  \ite In Section \re{sec:overfiiting}, we provide an interpretation of the bias in the optimal cleaning procedure (\ie of the lack of convergence of $s_k\ucleaned$ as estimator of $s_k$)  in terms of the overfitting factor of the estimator, out of $X$,  of a certain projection of $Y$.
  \ent


\subsection{Notations} Throughout this text, for $M$ a matrix, $M'$ denotes the transpose of $M$. For $Z$ a random variable,   $\E Z$   denotes the expectation of $Z$.  

 Here,  error terms in   approximations     depend on  the parameters $n$, $p$, $T$ and $\Si$ of the problem,  on the complex number $z$ and on the randomness. We will suppose that $\f nT$, $\f pT$, the operator norm of $\Si$ and $|z|$ are bounded by a constant $\mathfrak{M}$ and use the notation $$O\lf(\ff{T\, |\Im z|^k}\ri)=\f{O(1)}{T\, |\Im z|^k}$$ for error terms with the following definition: for  $Z=Z_{n,p,T,\Si,z}$ a complex  random variable  depending on $n,p,T,\Si,z$, we write $Z=O(1)$ if there exists $C$ depending only on  $\mathfrak{M}$ \st   $\E e^{|Z/C|^2}\le 2$, 
 \ie if $Z$ is Sub-Gaussian\footnote{Definition and basic properties of Sub-Gaussian variables can be found in \cite[Sec. 2.5]{Vershynin}.} with Sub-Gaussian norm controlled by  $\mathfrak{M}$.

   \section{Main results and Algorithms}\la{sec:modelresultsalgos}

 \subsection{Model}\la{subsec:modelresultsalgos}
 Let    $n\le p$ and  let  $(X,Y)\in \R^n\ti \R^p$ be  a  pair of  (column) random vectors\footnote{We suppose that $n\le p$ to avoid spurious null eigenvalues for the matrix $\cC\cC^*$.} such that $$\bbm X\\ 
Y\ebm\sim\cN(0, \Si)$$ for a given $\Si=\bpm\cA & \cC\\ \cC'&\cB\epm\in \R^{(n+p)\ti(n+p)}$ symmetric and non negative definite.
 
 We are  interested in the estimation  of   the true  cross-covariance  matrix $$\cC=\E XY'\in \R^{n\ti p}$$ out of its empirical version  \be\la{defbE}\bE:=\ff T\bX\bY' \in \R^{n\ti p},\ee where 
  \be\la{defXYFCM}\bX:=\bbm X(1)&\cd& X(T)\ebm\in \R^{n\ti T}\AND \bY:=\bbm Y(1)&\cd& Y(T)\ebm\in \R^{p\ti T}\ee are defined thanks to a sequence  \be\la{eq:realXY0}(X(1),Y(1)), \ld, (X(T),Y(T))\ee  of  independent copies of $(X,Y)$.

   More precisely, we are looking for a \emph{Rotationally Invariant Estimator}   $\mathbf{C}_{XY,\op{RIE}}$  of $\cC$, \ie an estimator 
   deduced from the estimator $\bE$ from \eqre{defbE} by \emph{changing} (following other papers on close questions, we say \emph{cleaning}) its singular values but not changing its singular vectors, so that for any $V,W$ orthogonal matrices,  if $\bX$ and $\bY$ are respectively changed into $V\bX$ and $W\bY$, then $\mathbf{C}_{XY,\op{RIE}}$ is changed into $V\mathbf{C}_{XY,\op{RIE}}W'$.

    Let us  introduce the SVD of $\bE$. 
    We  
    set  \be\la{eq1SVD}\bE\;= \; \sum_{k=1}^{n} s_k\bu_k\bv_k'\;= \; \underbrace{\bbm \bu_1 &\cd & \bu_n\ebm}_{:=\bU} 
    \diag(s_1, \ld, s_n) 
    \underbrace{\bbm \bv_1 &\cd & \bv_n\ebm'}_{:=\bV'}\ee
    for some $s_1, \ld, s_n\ge 0 $  
    and two orthonormal column vectors systems $\bu_1, \ld, \bu_{n}\in \R^n$,    
 and   $\bv_1, \ld, \bv_{n}\in \R^p$. 
  
  Thus our estimator will have the form $$\mathbf{C}_{XY,\op{RIE}}=\bU\diag(s\ucleaned_1, \ld, s\ucleaned_n)\bV'$$ 
  and the \emph{cleaned   singular values}  \be\la{Lancieux_221218}s\ucleaned_1, \ld, s\ucleaned_n\ee will be considered  \emph{optimal} when    solving the   optimization problem
  \be\la{opt_problem}\min_{s\uclean_1, \ld, s\uclean_n }\|\bU\diag(s\uclean_1, \ld, s\uclean_n)\bV'-\cC\|_{\text{F}} ,\ee
where  the Frobenius norm
  $\|\cdot\|_{\text{F}}$  has been defined at \eqre{opt_problem00}.  
   Let us introduce the (implicitly depending on $z\in \C\bck\R$) random variables  
  \begin{align}\la{Notation:Intro250700AA} &G=G(z):=\ff T\Tr  \bG,\qquad &&L=L(z):=\ff T \Tr \bG \bC_{XY}\cC' \end{align} 
  for $\bG$ (resp. $\wtG$, that we shall also use below) the resolvent, estimated at $z^2$,  of $\bE\bE'$ (resp. of $\bE'\bE$) defined through \be\la{eq:defGwtG}\bG :=\lf(z^2-\bE\bE'\ri)^{-1},\qquad \wtG :=\lf(z^2-\bE'\bE\ri)^{-1}.\ee  
  Note that $G(z)$ can be computed with the observed data, so that  by \eqre{0407184IO00AA} bellow, the function $L(z)$  is all one needs to compute the cleaned singular values $s\ucleaned_k$.   For this reason, it  is called the \emph{oracle function}.
  \beg{propo}\la{propo:sol_optim}The solutions of the optimization problem \eqre{opt_problem} satisfy  \be\la{0507181}s\ucleaned_k\;=\;\bu_k'\cC\bv_k\;\approx   \;\f{\Im L(z)}{\Im (zG(z))} \quad\text{ for $z=s_k+\ii \eta$, with $\eta\ll 1$}
  ,\ee where the functions  $L(z)$ and $G(z)$ are defined    at \eqre{Notation:Intro250700AA}. More precisely, 
  for any $\eps>0$ \st $[s_k-\eps, s_k+\eps]\cap\{s_1, \ld, s_n\}=\{s_k\}$, \be\la{0407184IO00AA}s\ucleaned_k\;=\;\bu_k'\cC\bv_k\;=\;\lim_{\eta\to 0}\f{\int_{s_k-\eps}^{s_k+\eps}\Im L(x+\ii\eta)\ud x}{\int_{s_k-\eps}^{s_k+\eps}\Im ((x+\ii\eta)G(x+\ii\eta))\ud x} .\ee
  \end{propo}

  \brem[Impact of a $o(1)$ error on $s\ucleaned_k$]\la{small_error_effect}  Equation \eqre{0407184IO00AA} provides us with an exact formula for $s\ucleaned_k$, that, once an explicit approximation for the function $L$   obtained, we shall convert into an approximate formula $s\ucleanedalgo_k$ for $s\ucleaned_k$ in \eqre{16081810h}. One can wonder what  the effect of this approximation is on the optimality (in Frobenius norm, with asymptotic   probability tending to one) of  estimator $\mathbf{C}_{XY,\op{RIE}}=\bU\diag(s\ucleaned_1, \ld, s\ucleaned_n)\bV'$. 
  We are thus interested in the matrix error term $$\bU\diag(s\ucleanedalgo_1, \ld, s\ucleanedalgo_n)\bV'-\bU\diag(s\ucleaned_1, \ld, s\ucleaned_n)\bV'$$ 
  The question is: \emph{How small must the errors  $s\ucleanedalgo_k-s\ucleaned_k$ be
  for the matrix error term above to  be negligible with respect to the true cross-covariance matrix $\cC$ for the Frobenius norm?}
  Given $\bU$ and $\bV$ are orthogonal matrices, the Frobenius norm of the matrix error term is
   $$ \lf(\sum_{k=1}^n (s\ucleanedalgo_k-s\ucleaned_k)^2 \ri)^{1/2}
  .$$
  Under the sole hypothesis that the operator norm of $\Si$ is bounded by the constant $\mathfrak{M}$, the Frobenius  norm of $\cC$ has order \emph{a priori} $\sqrt{n}$  so that  the error on $\mathbf{C}_{XY,\op{RIE}}$ is  negligible as soon as $$\sum_{k=1}^n (s\ucleanedalgo_k-s\ucleaned_k)^2=o(n).$$
Regimes where the actual Frobenius  norm of $\cC$ has lower order (e.g.  when $\cC$ is simply null) have to be the object of specific studies. 
  \erem

  \subsection{Estimations of the oracle function $L(z)$}The problem with   Formula \eqre{0507181} is that while the function $G(z)$ is explicit from the data $\bX,\bY$, the definition of the function $L(z)$ involves the unknown true cross-covariance matrix $\cC$.  In Theorems \re{th:main} and \re{th:main2}, we   give   asymptotic approximations of $L(z)$ that can be estimated from the data alone, as is the case of the Ledoit-P\'ech\'e estimator for covariance matrices \cite{RIE0}.

   Let us introduce the  random variables  \begin{align}\la{Notation:Intro2507} &H:=\ff T\Tr \bG\bE\bE',\quad A:=\ff T \Tr \bG \bC_X, \quad B:=\ff T \Tr \wtG \bC_Y, \quad \Tta :=z^2\f{AB}{1+H} 
   \end{align} for $\bG$,  $\wtG$  as in 
    \eqre{eq:defGwtG} and   $\bC_X,\bC_Y$ the empirical covariance matrices of $X$ and $Y$ defined by \be\la{eq:defCXCY}\bC_X:=\ff T \bX\bX',\qquad \bC_Y:=\ff T \bY\bY'.
\ee 

The following result makes the function $L$ of  \eqre{Notation:Intro250700AA} explicit from the data alone, allowing a practical implementation of   Formula \eqre{0507181} for the RIE. 
  
  \bTh[Oracle function estimation I]\la{th:main}
  The function  $L$ of  \eqre{Notation:Intro250700AA}   satisfies 
  \be\la{eq:LmadeExplicit} L(z)=\f{H(z)-\Tta(z)}{1+H(z)-\Tta(z)} +O\lf(\ff{T\, |\Im z|^5}\ri).\ee  \eTh

\brem[Case where $T\gg n,p$]\la{rmk:boundednp} In the case where, as $T$ tends to infinity, $n$ and $p$ stay bounded, it can easily be seen that $L\approx H$, so that $$\mathbf{C}_{XY,\op{RIE}}\approx \bE.$$ Indeed, the estimate  $L\approx H$ follows for example from the formulas (true  for large $|z|$):  \beq \f TnL&=&\sum_{k\ge 1}\f{z^{-2k}}n\Tr \lf((\bE\bE')^{k-1}\bE\cC'\ri) \\
   \f TnH&=&\sum_{k\ge 1}\f{z^{-2k}}n\Tr \lf((\bE\bE')^{k-1}\bE\bE'\ri)
   \eeq and from standard complex analysis. 
\erem

In the particular case where the covariance matrices of $X$ and $Y$ are both identity matrices,  $\mathbf{C}_{XY,\op{RIE}}$ is  in fact an estimator of the \emph{cross-correlation} matrix of $X$ and $Y$, and \eqre{0507181} can be simplified into \eqre{explicit130818final2}, a formula leading to an algorithm 
with lower computational complexity (see Remark \re{compared_complexities}).
For $z\in \C\bck\R$, let  \be\la{Notation:Intro2507Cor} 
K:=\lf(\f{p-n}T+z^2G\ri)G(1+H)^2. \ee

  \bTh[Oracle function estimation II]\la{th:main2} Suppose that the true covariance matrices   of $X$ and $Y$ are  respectively  $I_n$  and   $I_p$. Then,   the function  $L$ of  \eqre{Notation:Intro250700AA}   satisfies 
   \be\la{explicit130818final2}L(z)=\f{1+2H(z)-\sqrt{1+4K(z)}}{2(1+H(z))}+O\lf(\ff{T\, |\Im z|^5}\ri)\ee
 for $\sqrt{\cdote}$ the analytic version of the square root on $\C\bck(-\infty, 0]$ with value $1$ at $1$.\eTh
  
  \subsection{Algorithmic consequences}\la{sec:algos}
  Formula \eqre{0507181}
  gives an expression for the cleaned singular values of the cross-covariance matrix, \ie for the RIE of this matrix. The  function $G(z)$ is explicit from the data $\bX,\bY$, as well as the approximation of $L(z)$ given by formulas \eqre{eq:LmadeExplicit} and \eqre{explicit130818final2} above. Choosing $\eta=(npT)^{-1/12}$ for the "small $\eta$" (and using the formula  $H=z^2G-n/T$, for $H$ as in \eqre{Notation:Intro2507})  leads to the explicit implementation formula 
\be\la{16081810h}s\ucleanedalgo_k\;:=\;s_k\ti  \f{\Im L(z)}{\Im H(z)} \quad\text{ for $z=s_k+\ii (npT)^{-1/12}$}.
\ee

\brem\la{rem:actual_error_term} As explained in Remark \re{small_error_effect},  we need  the error $s\ucleanedalgo_k-
s\ucleaned_k$ (squared, and averaged over $k$) to be $o(1)$, so that with the $O\lf(\ff{T\, |\Im z|^5}\ri)$ error terms of 
 \eqre{eq:LmadeExplicit} and \eqre{explicit130818final2}, we should have slightly increased the imaginary part of $z$ in \eqre{16081810h}, to have $T\, |\Im z|^5$   large. It happens that  in practice, the algorithms below work well with $(npT)^{-1/12}$ as imaginary part of $z$. 
  In fact, we believe that,
   following the method developed by Erd\H{o}s, Yau and co-authors (see e.g. \cite{ErYauS,ErYau,FloAntti}), our local  laws in Theorems \re{th:main} and \re{th:main2} can be improved roughly up  to the scale $T^{-1}$, \ie that the error terms, in 
   \eqre{eq:LmadeExplicit} and \eqre{explicit130818final2}, are in fact controlled essentially by $(T \Im z)^{-1}$.  
   This conjecture has been tested in Section \re{sec:numsimoracle} (see the caption of \fig{fig1_num}).
  \erem

Using the approximation of $L(z)$ given by formula \eqre{eq:LmadeExplicit}, we get  the first  algorithm below, whose complexity is kept reasonable thanks to the following. With \be\la{def:SVDbE1908}\bE\;= \;   \bbm \bu_1 &\cd & \bu_n\ebm  \diag(s_1, \ld, s_n) \bbm \bv_1 &\cd & \bv_n\ebm'\ee
the SVD of $\bE$, where the orthonormal system $ \bv_1 ,\ld , \bv_n$ of $\R^p$ is completed to an orthonormal basis $ \bv_1 ,\ld , \bv_p$, 
we have \be\la{eq:algoFormulas1}H(z)=\ff T\sum_{\ell=1}^n \f{s_\ell^2}{z^2-s_\ell^2},\;
A(z)=\ff T\sum_{\ell=1}^n \f{a_{\ell}}{z^2-s_\ell^2},\; B(z)=\ff T\lf(\sum_{\ell=1}^n \f{b_{\ell}}{z^2-s_\ell^2}+z^{-2}b_{[n+1:p]}\ri)\ee
for  \be\la{eq:algoFormulas2}a_{\ell}:=\bu_\ell'\bC_X\bu_\ell,\qquad b_{\ell}:=\bv_\ell'\bC_Y\bv_\ell,\qquad b_{[n+1:p]}:=\sum_{\ell=n+1}^p\bv_\ell'\bC_Y\bv_\ell \ee
so that the functions  $H$, $A$, $B$ and $\Tta$ from \eqre{Notation:Intro2507} can be computed without any matrix inversion (nor any matrix product) once the SVD of $\bE$ has been computed, which has only to be done once in the algorithm.

\begin{algo}{Algorithm 1: Optimal cleaning for cross-covariance  matrices}
\emph{Input:} $\bX\in \R^{n\ti T}$, $\bY\in \R^{p\ti T}$ with $n\le p$.

\emph{Output:} cleaned singular values $s\ucleanedalgo_1, \ld, s\ucleanedalgo_n$.
\\
\beg{enumerate}
\ite Compute $\bE=\ff T\bX\bY'$, $\bC_X=\ff T\bX\bX'$, $\bC_Y=\ff T\bY\bY'$
\ite Compute the SVD 
of $\bE$
    \ite Compute the vectors $(a_{\ell})_{\ell=1, \ld, n}$ and  $(b_{\ell})_{\ell=1, \ld, n}$ and the number $b_{[n+1:p]}$ using   \eqre{eq:algoFormulas2}
    \ite For each $k\in \{1, \ld, n\}$,
    \bgt\ite  set $z=s_k+\ii (npT)^{-1/12}$  for $s_k$ the $k$-th singular value of $\bE$
    \ite compute $H$, $A$, $B$ using   \eqre{eq:algoFormulas1}  
    \ite  compute $\Tta=z^2\f{AB}{1+H}$ and $ L=1-\ff{1+H-\Tta} $
    \ite compute
   $$s\ucleanedalgo_k\;=\;s_k\ti  \Big(\f{\Im L}{\Im H} \Big)_+\quad\text{(with $x_+:=\max\{x,0\}$)}
$$\ent
\vspace{-.1cm}\ite possibly: apply the isotonic regression algorithm to the $s\ucleanedalgo_k$
\en{enumerate}
\end{algo}

One can also write an algorithm based on \eqre{explicit130818final2} instead of \eqre{eq:LmadeExplicit}, with slightly lower computational complexity, but only works when  the true covariance matrices of $X$ and $Y$ are both identities (which can be the case in practice, when the associated data has been made standard in a preprocessing):
 
\begin{algo}{Algorithm 2: Optimal cleaning for cross-correlation  of signals with identity covariance matrix}
\emph{Input:} singular values $s_1, \ld, s_n$ of $\bE=\ff T\bX\bY'$ for $\bX\in \R^{n\ti T}$, $\bY\in \R^{p\ti T}$.

\emph{Output:} cleaned singular values $s\ucleanedalgo_1, \ld, s\ucleanedalgo_n$.
\\

 For each $k\in \{1, \ld, n\}$,
   \beg{enumerate}\ite set $z=s_k+\ii (npT)^{-1/12}$
    \ite compute $$ G=\ff T \sum_{\ell=1}^n \ff{z^2-s_\ell^2},\quad H=z^2G-n/T,\quad K=\lf(\f{p-n}T+z^2G\ri)G(1+H)^2$$ and $$ L=\f{1+2H-\sqrt{1+4K}}{2(1+H)}$$
       \ite compute
   $$s\ucleanedalgo_k\;=\;s_k\ti    \Big(\f{\Im L}{\Im H} \Big)_+
$$ 
\vspace{-.1cm}\ite possibly: apply the isotonic regression algorithm to the $s\ucleaned_k$
\en{enumerate}\end{algo}

\brem Ledoit-P\'ech\'e's RIE  is not working when the ratio $q$ of the signal size by the sample size is too close to 1 because its formula involves a division by $q-1$. No such singularity appears here. 
\erem

\brem[Compared computational complexities of Algorithms 1 and 2]\la{compared_complexities} Using the classical linear algebra operations computational complexity estimates (multiplication, inversion and singular value decomposition of $O(N)\times O(N)$ matrices have $O(N^3)$ complexity), we see that both algorithms have complexities $O(T^3)$. That being said, Algorithm 2 involves less matrix multiplications than Algorithm 1, given it does not use the numbers defined at \eqre{eq:algoFormulas2}. It follows  that when both algorithms apply, Algorithm 2 needs less computation time, as Figure \re{fig_comparison_complexities} illustrates.\erem
 \begin{figure}[ht]
\centering
\includegraphics[scale=.6]{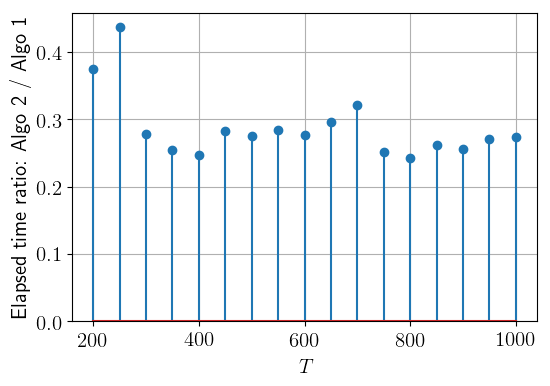}
\caption{\textbf{Compared computational complexities of Algorithms 1 and 2:} for each value of $T$, we compute the ratio of the total time elapsed to run both algorithms with 100 simulations of the model described in Section \ref{sec:numsimoracle}
 with fixed $n/T=0.4$ and $p/T=0.7$.}\label{fig_comparison_complexities}
\end{figure}

\subsection{Cleaned vs empirical vs true  singular values: some exact formulas}\la{sec:comparisonEmpCleanTrue} 
As explained above in Section \re{section_optimality}, due to the unavoidable error in the singular vectors, the fact that our estimator   realizes the optimal of \eqre{intro_optimal_eq} does not imply that the cleaned singular values should be close to the true ones.
  Precisely, we show in Proposition \re{prop_unbiased_estimator} 
  that on average, 
\be\la{rough_estimate_cleaned_true_emp2}s_k\ucleaned\;<\;s_k\utrue\;<\;s_k,\ee
for $s_k\utrue$ the true singular values (the left inequality being conditional to $s_k\utrue>0$).

An analogous phenomenon for rotationally invariant estimators of covariance matrices is explained   in \cite[Section 6.3]{BBPo}.

\beg{propo}\la{prop_unbiased_estimator}
Let $\lam_k\Xtrue$ (resp. $\lam_l\Ytrue$) denote the  eigenvalues of the true covariance matrix of $X$ (resp. $Y$) and let  $\bC_X$, $\bC_Y$ denote the empirical covariance matrices of $X$ and $Y$ from \eqre{eq:defCXCY}. Then the following equations hold:
\be\la{precise_estimate_cleaned_true_emp1}
\E  \sum_{k=1}^n s_ks_k\ucleaned\;=\;  \sum_{k=1}^n (s_k\utrue)^2,
\ee
\be\la{precise_estimate_cleaned_true_emp2}
\E  \sum_{k=1}^n s_k^2\;=\; \lf(1+\ff T\ri) \sum_{k=1}^n (s_k\utrue)^2+\f2T\sum_{k=1}^n\lam_k\Xtrue\sum_{l=1}^p\lam_l\Ytrue
\ee and
\be\la{precise_estimate_cleaned_true_emp3}
   \sum_{k=1}^n (s_k\utrue)^2\;=\;\ff{1+T^{-1}-2T^{-2}}\E\lf[ \sum_{k=1}^n s_k^2-\ff{T}\Tr \bC_X\Tr \bC_Y\ri].\ee 
\en{propo}
 
Note that the rough estimate of \eqre{rough_estimate_cleaned_true_emp2} follows from \eqre{precise_estimate_cleaned_true_emp1} and \eqre{precise_estimate_cleaned_true_emp2}: 
it follows from \eqre{precise_estimate_cleaned_true_emp2} that on average, 
$$s_k\utrue\,<\,s_k,$$ whereas it follows from \eqre{precise_estimate_cleaned_true_emp1} that on average, 
$$s_k\ti s_k\ucleaned\;\approx\;(s_k\utrue)^2.$$

\subsection{Interpretation of the cleaning in terms of overfitting}\la{sec:overfiiting}\emph{Overfitting} is a very common issue in  machine learning. It refers to the problem that any model is fitted, \emph{in sample}, on noisy data, which can degrade its \emph{out of sample} performance if the fit has been significantly impacted by the random specificity of the noise. In this section, we relate the cleaning procedure to this problem, proving that in some contexts, 
$$\frac{\text{Out-of-sample-performance}}{\text{In-sample-performance}}\approx \frac{s_k\ucleaned}{s_k}.$$

Suppose to be given a sample $$(X(1),Y(1)),\ldots,(X(T),Y(T))$$ of observations  of a pair $(X,Y)\in\R^{p}\ti\R^n$ of vectors, where $X$ is a collection of factors thanks to which we want to explain   $Y$. 

Given this set of observations, if we observe an ``out-of-sample'' (oos) realization $X\oos$ of the factors,  a natural predictor\footnote{This predictor, a matched-filter, corresponds to   the normalized Ridge predictor with large $\lam$ (namely the large $\lam$ limit of  \cite[Eq. (3.47)]{Tibshi} times $\lam$).  We could also consider the OLS predictor, but notations are lighter this way.} of the corresponding $Y\oos$ is given  by 
$$ \op{Pred}_Y(X\oos):=\bERF X\oos=\sum_k(s_k\ti X\oos'\bv_k)\bu_k,$$ where $ \sum_{k=1}^n s_k\bu_k\bv_k'$ is the SVD of the in-sample cross-covariance matrix $$\bERF:=\ff T\sum_t Y(t)X(t)'.$$

Each term of the previous sum defines a partial predictor\footnote{Note that we do not name this an \emph{estimator} but a \emph{predictor}: the purpose of $\op{Pred}_Y^{(k)}$ is not to approximate $Y$ the better but to exhibit  a positive  alignment  with  $Y$. This kind of object is widely used in finance (see Remark \re{rmk:financialinteroretation}), where for various reasons (risk aversion, volatility, liquidity, crowding, better predictability), an investor mights want to focus on some specific directions in the market.} $$\op{Pred}_Y^{(k)}(X\oos):=(s_k\ti X\oos'\bv_k)\bu_k .$$

Let us now focus on the  overlap of these predictors with the true values of $Y$.

\textbf{Out of  sample overlap:} it is given  by 
\beq Y\oos\cdot \op{Pred}_Y^{(k)}(X\oos)&=&s_k(X\oos'\bv_k)(\bu_k'R\oos)\\ &=&s_k\bu_k'Y\oos X\oos'\bv_k.\eeq
Over the out-of-sample time series  $$(X\oos(1),Y\oos(1)),\ldots,(X\oos(T\oos),Y\oos(T\oos)),$$
it   averages out to the \emph{mean out of sample overlap} $\meanoosovl$ given by 
\beq \meanoosovl&:=&\ff{T\oos}\sum_{t=1}^{T\oos}s_k\bu_k'Y\oos(t) X\oos(t)'\bv_k\\ &=&
s_k\bu_k'\bERFoos\bv_k
\eeq for $$\bERFoos:= \ff{T\oos}\sum_{t=1}^{T\oos} Y\oos(t) X\oos(t)'.$$
By \eqre{def:sikskIntro} and the concentration of measure  Lemma \re{lemma_concentration_overfitting},    we get
\be\la{osoverlap1218}\meanoosovl= s_ks_k\ucleaned +O\lf(\ff{\sqrt{T\oos}}\ri).\ee

\textbf{In sample overlap:} it is given, at each date $t$ of the sample, by 
$$Y(t)\cdot \op{Pred}_Y^{(k)}(X(t))=s_k(X(t)'\bv_k)(\bu_k'Y(t))=s_k\bu_k'Y(t)X(t)'\bv_k,$$
which, by  \eqre{def:sikskIntro},  averages out,  in   sample,   to the \emph{mean in sample overlap} $\meanisovl$ given by  \be\la{isoverlap1218}\meanisovl:=\ff T\sum_{t=1}^Ts_k\bu_k'Y(t)X(t)'\bv_k=s_k\bu_k'\bERF \bv_k=s_k^2.\ee

\textbf{Out of  sample / in sample:} From \eqre{osoverlap1218} and \eqre{isoverlap1218}, we deduce the following nice relation between the overfitting  and the cleaning, illustrated at \fig{Fig_Overfitting_Factor}: 
\be\la{osoisoverlap1218}\frac{\meanoosovl}{\meanisovl}\approx \frac{s_k\ucleaned}{s_k}.\ee
         \newcommand{\scaleOFTFact}{.7}
    \begin{figure}[ht]
\centering
\includegraphics[scale=\scaleOFTFact]{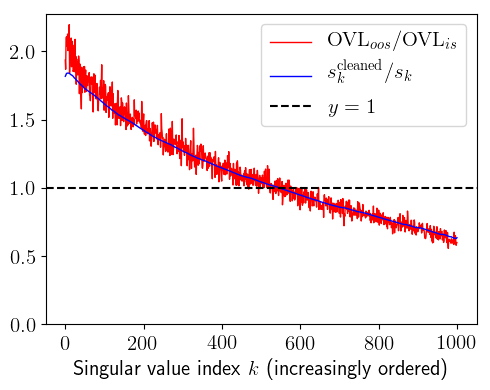} 
\caption{\textbf{Overfitting and cleaning}. Illustration of \eqre{osoisoverlap1218}: LHT vs RHT. Here, $Y=AX+\text{noise}$, with $X$ standard Gaussian vector, $A\in \R^{n\ti p}$ fixed, $n=p=1000$, $T=10000$, $T\oos=1000$.}\label{Fig_Overfitting_Factor}
\end{figure}

\brem If, instead of considering the  partial predictor $\op{Pred}_Y^{(k)}(X\oos)$, we consider sums, over $k$, of such predictors, then the previous ratio can still be expressed thanks to the numbers  $s_k\ucleaned$ and $s_k$.
\erem

 \brem[Investment strategies interpretation]\la{rmk:financialinteroretation} In the case where $Y$ is the vector of returns of a collection of financial assets and $X$   a collection of factors we want to build an investment strategy on, 
the vectors $$\pi_k:=\op{Pred}_Y^{(k)}(X\oos)$$ (as well as linear combinations of such vectors) correspond to portfolios constructed thanks to the factors $X$: each component of $\pi_k$ is the  (positive or negative)  amount of money invested in the corresponding asset. In this context,  
   the out-of-sample overlaps 
   $$Y\oos\cdot \pi_k$$ are simply be the \emph{realized gains} of these strategies, whereas the mean in-sample overlaps are the \emph{predicted gains} of these strategies. 
   
 Of course, realized gains are usually different, even on average, from predicted gains. This phenomenon can be seen as a consequence of \emph{overfitting} (or \emph{in-sample bias}). For  the simple model presented in this section,  
 \eqre{osoisoverlap1218}  relates their ratio  to the singular values cleaning procedure via the formula:
 $$\frac{\text{realized gains}(\pi_k)}{\text{predicted gains}(\pi_k)}\approx\frac{s_k\ucleaned}{s_k}.$$
\erem 

 \section{Numerical simulations}\la{sec:numsim}

  \subsection{Oracle  estimation}\la{sec:numsimoracle}The cornerstones of this work are 
     \eqre{eq:LmadeExplicit} and \eqre{explicit130818final2} from Theorems \re{th:main} and \re{th:main2}: these formulas allow us to approximate the (unknown) oracle function $L(z)$ by some functions that are explicit from the data.
We conducted numerical simulations to verify these formulas  for various models (\ie various choices of $\Si$), all  confirming their accuracy. In \fig{fig1_num}, we present the \emph{relative differences} 
\be\la{RefDif1}\f{\lf| L(z)-\f{H(z)-\Tta(z)}{1+H(z)-\Tta(z)}\ri|}{| L(z)|}\ee 
and 
\be\la{RefDif2}\f{\lf| L(z)-\f{1+2H(z)-\sqrt{1+4K(z)}}{2(1+H(z))}\ri|}{| L(z)|}\ee 
for 
  \be\la{RefDifModel}\bpm X\\ Y\epm\sim \cN(0,\Si)\quad\text{ with }\quad\Si=\begin{pmatrix} I_n&\cC\\ \cC'&I_p\end{pmatrix}\ee
with $\cC$ a matrix with singular values  distributed according to the bi-modal density from the right graph in \fig{Fig_CFM_talk} and independent, Haar-distributed, left and right singular vectors, independent from the singular values. We see that both approximations of $L(z)$ are   very efficient and that  the approximation of $L(z)$ given at \eqre{eq:LmadeExplicit} is slightly better, which is confirmed by other simulations.

 \begin{figure}[ht]
\centering
\includegraphics[scale=.6]{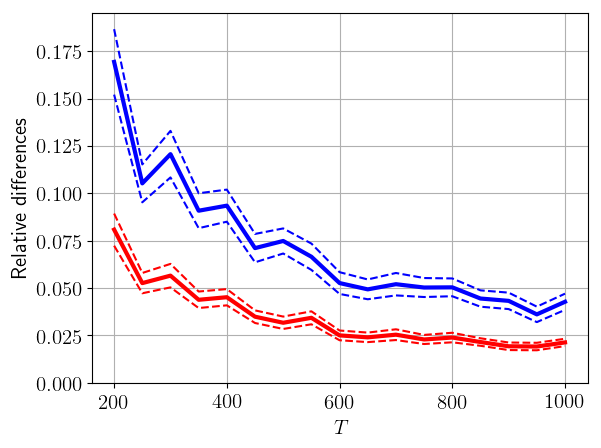}
\caption{\textbf{Validity of our oracle  estimations:} mean (out of 100 simulations, with 95\% confidence interval given by dashed lines) relative differences from \eqre{RefDif1} (in red) and \eqre{RefDif2} (in blue) for the model of \eqre{RefDifModel}, for various values of $T$ (in abscissa) with fixed $n/T=0.4$ and $p/T=0.7$ and for  $z=0.5+\ii(npT)^{-1/12}$. On the same simulations for the same value of $z$, we inferred  the exponent $q$ such that the error terms $ L(z)-\f{H(z)-\Tta(z)}{1+H(z)-\Tta(z)}$ and $L(z)-\f{1+2H(z)-\sqrt{1+4K(z)}}{2(1+H(z))}$ behave as $T^{-q}$ (note that the error terms plotted here  are the magnitudes of these error terms divided by $L(z)$). The inferred values are both $0.74$, which matches remarkably well the conjecture of Remark \re{rem:actual_error_term} (which  gives $q=3/4$).}\label{fig1_num}
\end{figure}

  \subsection{Effect of   cleaning}\la{sec:effect_of_the_cleaning}In \fig{Fig_CFM_talk_with_RIE} and \fig{Fig_scatter_plot}, we show the effect of   cleaning in the simulations   from  \fig{Fig_CFM_talk}.  
  In the  two graphs of \fig{Fig_scatter_plot}, we observe  that   for most values of $k$ (all but the smallest ones), we have  $s_k\ucleaned< s_k\utrue < s_k$, 
  as stated informally  in \eqre{rough_estimate_cleaned_true_emp2}. The  histograms of the right graph of \fig{Fig_CFM_talk_with_RIE} show the same, with less precision for the set of $k$'s for which this is true.
   \newcommand{\scaleCCvsRIE}{.5}
\begin{figure}[ht]
\centering
\includegraphics[scale=\scaleCCvsRIE]{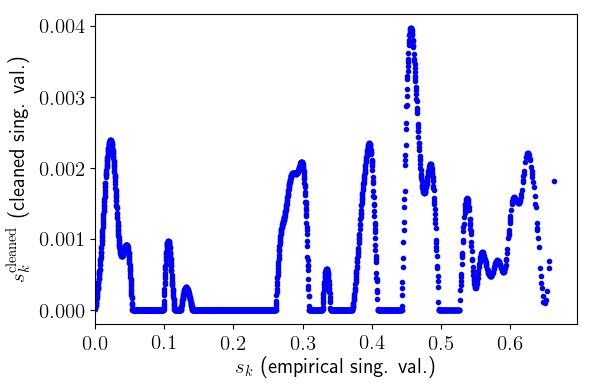}
\includegraphics[scale=\scaleCCvsRIE]{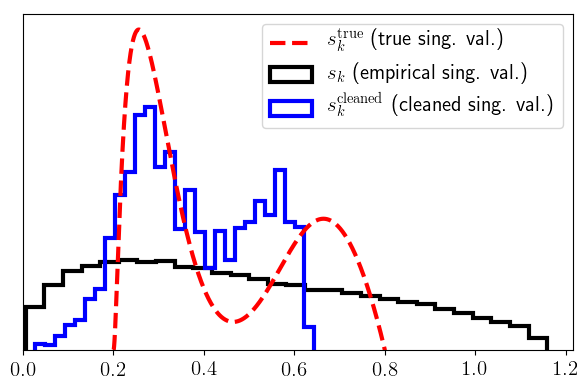}
\caption{\textbf{Cleaned vs empirical singular values} for the simulation from \fig{Fig_CFM_talk}. \textbf{Left:}  as one could expect from a good estimator, in the null model, most of the singular values are turned to approximately $0$. 
\textbf{Right:} same as in \fig{Fig_CFM_talk},  with the cleaned singular values histogram added. The lack of monotonicity in the left graph is the reason why we added the  isotonic regression as optional last step in our algorithm.}\label{Fig_CFM_talk_with_RIE}
\end{figure}
\begin{figure}[ht]
\centering
\includegraphics[scale=\scaleCCvsRIE]{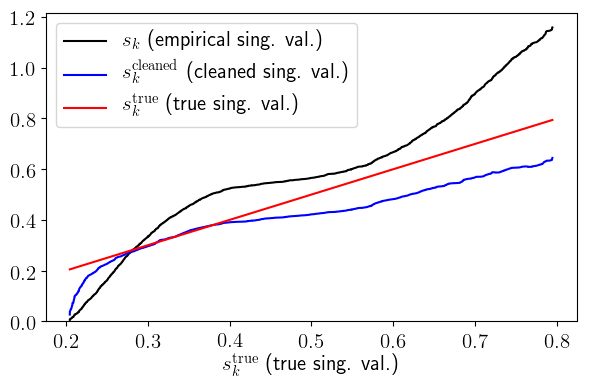}
\includegraphics[scale=\scaleCCvsRIE]{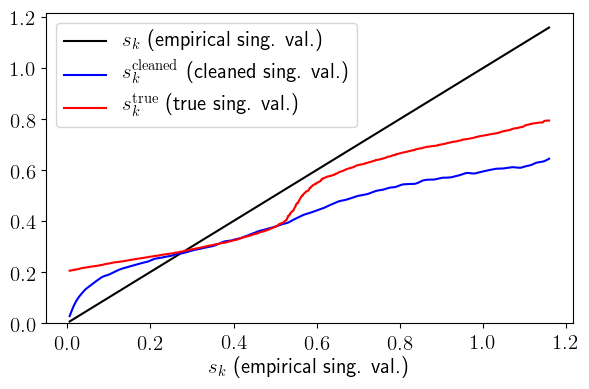}
\caption{Two alternative viewpoints on the right graph of \fig{Fig_CFM_talk_with_RIE}: empirical, cleaned and true singular values as functions of true (left) or of empirical  (right) singular values.}\label{Fig_scatter_plot}
\end{figure}

  \subsection{Compared performance with   empirical and Ledoit-P\'ech\'e's estimators}
 \subsubsection{Numerical simulations}
  We have implemented Algorithms 1 and 2  from the present paper\footnote{Both  give approximately the same result when $X$ and $Y$   have identity covariance matrices, 
  so we shall focus on Algorithm 1 in this section.} for various models, \ie various choices of the true total covariance matrix $\Si$ 
  \st
  $$\bpm X\\ Y\epm\sim \cN(0,\Si).$$ We then compared their performance to that of 
   the empirical  estimator  $\bE$, thanks to the relative distances  
 \be\la{def:distances_quotient}\text{RD}^{\text{Algo}/\text{Emp}}_{\text{F}}=\f{\lf\|\whcC\uAlgo-\cC 
 \ri\|_{\text{F}}}{\lf\|\bE-\cC \ri\|_{\text{F}}},\qquad \text{RD}^{\text{Algo}/\text{Emp}}_{\text{op}}=\f{\lf\|\whcC\uAlgo-\cC 
 \ri\|_{\text{op}}}{\lf\|\bE-\cC \ri\|_{\text{op}}},\ee where $\whcC\uAlgo$ denotes the estimator of $\cC$ obtained with Algorithm 1 and where $\|\cdot\|_{\text{op}}$ denotes the operator norm\footnote{Our algorithm is optimal, in the sense given in this paper, for the Frobenius norm, but the operator norm of the error is of course also interesting.}.
 We also 
compared  with   the $n\ti p$ upper-right  corner $\whcC\uLedoitPeche$ of Ledoit-P\'ech\'e's estimator of the total covariance matrix $\Si$ 
(the estimator from \cite{RIE0}, which assumes $O(n+p)$-invariance), thanks to the relative distances 
  \be\la{def:distances_quotientLP}
  \text{RD}^{\text{Algo}/\text{Ledoit-Peche}}_{\text{F}}=\f{\lf\|\whcC\uAlgo-\cC 
 \ri\|_{\text{F}}}{\lf\|\whcC\uLedoitPeche-\cC \ri\|_{\text{F}}}, \;   \text{RD}^{\text{Algo}/\text{Ledoit-Peche}}_{\text{op}}=\f{\lf\|\whcC\uAlgo-\cC 
 \ri\|_{\text{op}}}{\lf\|\whcC\uLedoitPeche-\cC \ri\|_{\text{op}}}
 \ee
 The  values of the quotients from \eqre{def:distances_quotient} and  \eqre{def:distances_quotientLP}    are  reported in Table \re{ref:my_table} and \re{ref:my_table_op} (and discussed  
in Section \re{sec:conc_persp} 
  below) 
  for the following models, all with  for $T=1000$, $n/T=0.4$ and $p/T=0.7$:
 \bgt
\ite Models (1) to (5):
  $$X\sim\cN(0,I_n) 
\;,\qquad Y=\cC'X+\si N,\; \text{ with $N\sim\cN(0,I_p)$ independent of $X$},$$
 where $\si^2=0.5$ and 
     $\cC$ has:
   \bgt\ite  independent Haar-distributed left and right singular vectors,      \ite  0\%, 10\%, 20\%, 30\% or 40\% (for respectively Model (1)\footnote{Model (1) corresponds in fact the null case from Figures \re{Fig_CFM_talk} and \re{Fig_CFM_talk_with_RIE}  even if $\si\ne 1$, given it can easily be seen that if $X$ or $Y$ is multiplied by a positive constant, then the outputs $s_k\ucleaned$ from our algorithms are also multiplied by this constant.},\ldots, Model (5)) of non zero singular values, distributed uniformly in $[0.2,0.5]$ (and independent of the singular vectors),\\
\ent
so that $\Si=\bpm I_n& \cC\\ \cC'&\cC'\cC+\si^2I_p\epm$.\\
 \ite Models (6) to (10): $\Si=HH'/(2m)$ for $H$ an $m\ti 2m$ matrix with i.i.d.  entries with common law $\mu$, which is either the standard Gaussian law (Model (6)) or  a symmetric heavy-tailed distribution with exponent $\al$ (specifically the symmetric law $\mu_\al$ such that   $\mu_\al(\R\bck[-x,x])=x^{-\al}$ for any $x\ge 1$), with    $\al=5$, $2.5$, $1.5$, $0.5$ for Models (7), (8), (9), (10) respectively. \beg{rmk}In Tables \re{ref:my_table} and \re{ref:my_table_op}, Models (6) and (7) always exhibit very close statistics: in our simulations, we observed that this phenomenon generalizes to any law $\mu$ with at least four moments. The same is true, with the same generalization, for models (11) and (12).
 \en{rmk}
  \ite For Models (11), $\ldots$, (15), use the $(n+p)\times(n+p)$ block decomposition $$\Si=\bpm \cA&\cC\\ \cC'&\cB\epm$$ of the matrix $\Si$ of respectively  Model (6), $\ldots$, (10)  and replace the covariance matrices $\cA$ and $\cB$ of $X$ and $Y$  by respectively $s_{\max}(\cC)I_n$ and $s_{\max}(\cC) I_p$, where $s_{\max}(\cC)$ denotes the largest singular value of $\cC$.
 \ent

 \begin{table}[ht]
  \vspace{.3cm}
 \begin{tabular}{|l|c|c|c|c|c|}
\hline Model & (1) & (2) & (3) & (4) & (5) \\
\hline 
Algo/Empirical 
&\FrameTable{0.01}{0.00004}
&\FrameTable{0.19}{0.00001}
&\FrameTable{0.26}{0.00001}
&\FrameTable{0.31}{0.00001}
&\FrameTable{0.35}{0.00001}\\
\hline 
Algo/Ledoit-P\'ech\'e
&\FrameTable{0.02}{0.00013}
&\FrameTable{0.53}{0.00003}
&\FrameTable{0.67}{0.00003}
&\FrameTable{0.75}{0.00003}
&\FrameTable{0.80}{0.00002}\\
\hline
\end{tabular}\vspace{.2cm}\\
\begin{tabular}{|l|c|c|c|c|c|}
\hline Model & (6) & (7) & (8) & (9) & (10) \\
\hline 
Algo/Empirical 
&\FrameTable{0.56}{0.00003}
&\FrameTable{0.56}{0.00003}
&\FrameTable{0.57}{0.00031}
&\FrameTable{0.47}{0.00682}
&\FrameTable{0.01}{0.00152}\\
\hline 
Algo/Ledoit-P\'ech\'e
&\FrameTable{0.95}{0.00002}
&\FrameTable{0.94}{0.00004}
&\FrameTable{0.85}{0.00058}
&\FrameTable{0.35}{0.00533}
&\FrameTable{0.01}{0.00152}
\\
\hline
\end{tabular}\vspace{.2cm}\\
\begin{tabular}{|l|c|c|c|c|c|}
\hline Model & (11) & (12) & (13) & (14) & (15) \\
\hline 
Algo/Empirical 
&\FrameTable{0.53}{0.00006}
&\FrameTable{0.53}{0.00006}
&\FrameTable{0.32}{0.00180}
&\FrameTable{0.10}{0.00062}
&\FrameTable{0.07}{0.00018}\\
\hline 
Algo/Ledoit-P\'ech\'e
&\FrameTable{0.97}{0.00008}
&\FrameTable{0.96}{0.00009}
&\FrameTable{0.54}{0.00392}
&\FrameTable{0.01}{0.00017}
&\FrameTable{0.00}{0.00003}
\\
\hline
\end{tabular}\vspace{.2cm}
\caption{\textbf{Frobenius norm comparisons:} confidence intervals for the means of the  ratios $\text{RD}^{\text{Algo}/\text{Emp}}_{\text{F}}$ of  \eqre{def:distances_quotient} (first row) and $\text{RD}^{\text{Algo}/\text{Ledoit-Peche}}_{\text{F}}$ from \eqre{def:distances_quotientLP} (second row) out of $10^4$ simulations.}\la{ref:my_table}
\end{table}

  \begin{table}[ht]
  \vspace{.3cm}
 \begin{tabular}{|l|c|c|c|c|c|}
\hline Model & (1) & (2) & (3) & (4) & (5) \\
\hline 
Algo/Empirical 
&\FrameTable{0.01}{0.00009}
&\FrameTable{0.39}{0.00009}
&\FrameTable{0.38}{0.00009}
&\FrameTable{0.37}{0.00008}
&\FrameTable{0.36}{0.00009}\\
\hline 
Algo/Ledoit-P\'ech\'e
&\FrameTable{0.04}{0.00028}
&\FrameTable{0.96}{0.00022}
&\FrameTable{0.92}{0.00017}
&\FrameTable{0.88}{0.00017}
&\FrameTable{0.86}{0.00014}\\
\hline
\end{tabular}\vspace{.2cm}\\
\begin{tabular}{|l|c|c|c|c|c|}
\hline Model & (6) & (7) & (8) & (9) & (10) \\
\hline 
Algo/Empirical 
&\FrameTable{0.46}{0.00014}
&\FrameTable{0.46}{0.00014}
&\FrameTable{0.58}{0.00217}
&\FrameTable{0.50}{0.01145}
&\FrameTable{0.02}{0.00185}\\
\hline 
Algo/Ledoit-P\'ech\'e
&\FrameTable{0.97}{0.00012}
&\FrameTable{0.95}{0.00035}
&\FrameTable{0.75}{0.00247}
&\FrameTable{0.43}{0.01058}
&\FrameTable{0.02}{0.00188}
\\
\hline
\end{tabular}\vspace{.2cm}\\
\begin{tabular}{|l|c|c|c|c|c|}
\hline Model & (11) & (12) & (13) & (14) & (15) \\
\hline 
Algo/Empirical 
&\FrameTable{0.45}{0.00011}
&\FrameTable{0.45}{0.00011}
&\FrameTable{0.52}{0.00050}
&\FrameTable{0.58}{0.00010}
&\FrameTable{0.58}{0.00010}\\
\hline 
Algo/Ledoit-P\'ech\'e
&\FrameTable{0.97}{0.00017}
&\FrameTable{0.97}{0.00017}
&\FrameTable{0.53}{0.00355}
&\FrameTable{0.01}{0.00013}
&\FrameTable{0.00}{0.00005}
\\
\hline
\end{tabular}\vspace{.2cm}
\caption{\textbf{Operator norm comparisons:} confidence intervals for the means of   the ratios $\text{RD}^{\text{Algo}/\text{Emp}}_{\text{op}}$ of  \eqre{def:distances_quotient} (first row) and $\text{RD}^{\text{Algo}/\text{Ledoit-Peche}}_{\text{op}}$ from \eqre{def:distances_quotientLP} (second row) out of $10^4$ simulations.}\la{ref:my_table_op}
\end{table}

\subsubsection{Comments}\la{sec:conc_persp}
On the examples from Tables \re{ref:my_table} and \re{ref:my_table_op}, 
our estimator always outperforms the empirical estimator by far, for the Frobenius norm as well as the operator norm.

For most of these examples, our estimator also outperforms significantly the upper-right corner of Ledoit-P\'ech\'e's estimator for both norms.   

Two factors seem to increase the advantage of our estimator over the two other ones considered here:
 \bgt\ite  Bayesian models with   prior distributions of the true total covariance matrix $\Si$   invariant under the action of $O(n)\ti O(p)$ defined by \be\la{1301181}(U,V)\cdot \Si=\bpm U&0\\ 0&V\epm \Si \bpm U'&0\\ 0&V'\epm,\ee (\ie Models (1) to (5)) are better fitted for our estimator than others, even the $O(n+p)$-invariant Bayesian models (case of Model (6)).\\
 \ite The sparser the true cross-covariance of a model, the higher our advantage over the two other algorithms: Models (7), (8), (9), (10) (or (12), (13), (14), (15)) are increasingly sparse (in the sense that a small part of the entries of $\cC$ contains most of its total mass).  Also,  the set of singular values of $\cC$ for  Models (1) to (5) are decreasingly sparse.
    \ent

   \beg{rmk} The $O(n)\ti O(p)$-invariance from \eqre{1301181} implies that the singular vectors of $\cC$ are Haar-distributed, but does not  imply that the singular vectors of $\cC$ are independent from the other observables (e.g. the eigenvectors of $\cA$ and $\cB$), hence does not define  Bayesian models where the right way to estimate $\cC$ is necessarily  rotationally invariant\footnote{Bayesian models where the right way to estimate $\cC$ is necessarily  rotationally invariant are those with prior distribution on $\Si$ invariant under the action of $O(n)^2\ti O(p)^2$ defined by $(U,W,V,K)\cdot \bpm \cA&\cC\\ \cC'&\cB\epm = \bpm U\cA U'&W\cC K'\\ K\cC'W'&V\cB V'\epm $.}. This means that for  Models (1) to (5), our estimator could be sub-optimal, and a cleaning of the singular vectors, based e.g. on the observation of the eigenvectors of $\cA$ and $\cB$, should possibly  also be performed. \en{rmk}

  \section{Proofs}\la{sec:proofs}
  \subsection{Proof of Proposition \re{propo:sol_optim}}\la{proof_of_prop_2p1}
 Proposition \re{propo:sol_optim} follows directly from both following claims. 
 
 \textbf{Claim 1.} The solution of the optimization problem \eqre{opt_problem} is given by \be\la{def:siksk}s\uclean_k=(\bU'\cC\bV)_{kk}=\bu_k'\cC\bv_k\;\quad  \text{ for } k=1, \ld, n.\ee
 
 \textbf{Claim 2.} For any $k=1, \ld, n$, $$\bu_k'\cC\bv_k\;=\;\lim_{\eta\to 0}\f{\int_{s_k-\eps}^{s_k+\eps}\Im L(x+\ii\eta)\ud x}{\int_{s_k-\eps}^{s_k+\eps}\Im ((x+\ii\eta)G(x+\ii\eta))\ud x}.$$
 
 \emph{Proof of Claim 1.}  Let $\tilde{\bV}$ be a $p\ti p$ orthogonal matrix with the same $n$ first columns as the $p\ti n$ matrix with orthogonal columns $\bV$. Then $$\bU\diag(s\uclean_1, \ld, s\uclean_n)\bV'=\bU\bbm\diag(s\uclean_1, \ld, s\uclean_n)&0_{n,p-n}\ebm \tilde\bV'$$ and,
  given the Frobenius norm is invariant by left and right multiplication by orthogonal matrices, the optimization problem \eqre{opt_problem} rewrites   
   $$\min_{s\uclean_1, \ld, s\uclean_n\ge 0}\|\bbm\diag(s\uclean_1, \ld, s\uclean_n)&0_{n,p-n}\ebm-\bU'\cC\tilde\bV\|_{\text{F}},$$
   \ie    \be\la{611211opti}\min_{s\uclean_1, \ld, s\uclean_n\ge 0}\|\diag(s\uclean_1, \ld, s\uclean_n)-\bU'\cC\bV\|_{\text{F}}.\ee As the squared  Frobenius norm of a matrix is simply the sum its squared entries, the 
    solution of \eqre{611211opti} is 
   given by the diagonal entries of $\bU'\cC\bV$, \ie by \eqre{def:siksk} . 
   \beg{rmk} The 
   two keys  to prove   Claim 1, 
   first 
   the 
   invariance 
   of 
   the 
   Frobenius norm under the left and right actions of the orthogonal group and second the fact that for any matrix $M$,  the diagonal matrix the closest to $M$ for the Frobenius norm is the diagonal matrix with the same diagonal entries as $M$, are, together,  specific to the Frobenius norm (at least among classical matrix norms). This is the reason why extending our results to other classical norms, such as the operator norm, 
   has so far remained  out of reach. That being said, simulations (see Table  \re{ref:my_table_op}) show that though possibly not optimal for the operator norm, using our estimator also makes sense (at least when compared to the empirical estimator) when the error is measured with this norm.\en{rmk}
   
    \emph{Proof of Claim 2.} The   $s\uclean_k$ of \eqre{def:siksk}  can be expressed as the Radon-Nikodym derivative  \be\la{0407183}s\uclean_k=\f{\ud \mE}{\ud \nu_\bE}(s_k),\ee
   for $\mE$  the null mass signed measure \be\la{0407181}\mE:=\ff{2n}\sum_{k=1}^n\bu_k'\cC\bv_k\lf(\del_{s_k}-\del_{-s_k}\ri)\ee and $\nu_\bE$  the symetrized empirical singular values distribution of $\bE$, defined by \be\la{0407182}\nu_{\bE}:=\ff{2n}\sum_{k=1}^n\lf(\del_{s_k}+\del_{-s_k}\ri).\ee

Equation \eqre{0407183} allows us, by \eqre{muoutofg}, to express $s\uclean_k$ thanks to the formula, true for any $\eps>0$ \st $[s_k-\eps, s_k+\eps]\cap\{s_1, \ld, s_n\}=\{s_k\}$, 
\be\la{0407184IO}s\uclean_k=\lim_{\eta\to 0}\f{\int_{s_k-\eps}^{s_k+\eps}\Im (g_{\mE}(x+\ii\eta))\ud x}{\int_{s_k-\eps}^{s_k+\eps}\Im (g_{\nu_\bE} (x+\ii\eta))\ud x} .\ee 
   By \eqre{eq:24aout20211}, 
   \beqy\nober g_{\mE}(z)&=&   \ff{n} \sum_{k=1}^n \f{s_k}{z^2-s_k^2}\bu_k'\cC\bv_k\\
 \nober     &=&\ff n \sum_{k=1}^n \f{s_k}{z^2-s_k^2}\Tr \bv_k\bu_k'\cC\\
  \nober          &=&\ff n \sum_{k=1}^n \f{s_k}{z^2-s_k^2}\Tr \cC'\bu_k\bv_k'\\
   \nober            &=&\ff n\Tr \lf(\cC' \sum_{k=1}^n \f{s_k}{z^2-s_k^2}\bu_k\bv_k'\ri)\\
   \nober           &=&\ff n\Tr \lf(\cC'\lf(z^2-\bE\bE'\ri)^{-1}\bE\ri)\\
   \la{1508181}           &=&\ff n\Tr \bG\bE\cC'
   \eeqy
      Similarly, by    \eqre{eq:24aout20212}, for $G$ as in \eqre{Notation:Intro250700AA},  
   \be\la{67182} g_{\nu_{\bE}}(z)=\f{T}nz G
    \ee
    Then, Claim 2 follows from \eqre{0407184IO}, \eqre{1508181}   and \eqre{67182}.

    \subsection{Proof of Theorem \re{th:main}}
 Let us introduce the   implicitly depending on $z\in \C\bck\R$ random variables  \be\la{Notation:Intro2507wtAwtB} \wtA:=\ff T\Tr \bG\cA\,\; \;\qquad 
  \wtB:=\ff T\Tr \wtG\cB.\ee

     Set
   $$g:=\E G,\;\; \;\; h:=\E H, \;\; \;\;\ell:=\E L,\;\; \;\; a:=\E A,\;\; \;\; \wta:=\E \wtA,\; \; \;\; b:=\E B,\;\;  \; \;\wtb:=\E \wtB
   .$$
  
  The following concentration of measure  lemma can   be proved using the Log-Sobolev inequality  satisfied by the standard Gaussian law (a detailed proof is given in Section \re{sec:measure_concentration}).
 \blem\la{lem:SubGaussianMainVariables} There is a   constant $c>0$, depending only on the bound $\mathfrak{M}$ of the hypothesis,    \st for any $z \in \C\bck\R$, we have, for any $t>0$, \beq \p\lf(\lf|G-g\ri|\ge t\ri)&\le&2\me^{-c(tT(\Im z)^4)^2}.
 \eeq 
 In other words, $G-g$ is a Sub-Gaussian random variable, with  Sub-Gaussian norm $$O\lf(\ff{T(\Im z )^4}\ri).$$
  
   Besides, the same is true for any of the   random variables $H-h$, $A-a$, $B-b$, $L-\ell$, $\Tta-\E\Tta $, $K- \E K$.
 \elem

   By this lemma, using the decomposition $$L=\ell +(L-\ell),$$  it suffices to prove that 
  \be\la{eq:1508188h} \ell=\f{h-\tta}{1+h-\tta} +O\lf(\ff{T\, (\Im z)^4}\ri).\ee  
   
   Then, the key of the proof is the following proposition, whose proof, based on the multidimensional  Stein formula for Gaussian vectors,  is postponed to Section \re{proofpropo:key_relations}.
      \beg{propo}\la{propo:key_relations}
We have \begin{align}\la{24071816h02} &h\;=\; \ell +\f{z^2  }{2}\lf( a\wtb+ b\wta  \ri)+     h\ell+O\lf(\ff{T\, (\Im z)^2}\ri)\\
  & \la{ODODOD}a(1-  \ell)\;=\;\wta(1+h)+O\lf(\ff{T\, (\Im z)^2}\ri) \\
  & \la{ODODODbis}b(1-  \ell)\;=\;\wtb(1+h)+O\lf(\ff{T\, (\Im z)^2}\ri)
   \end{align}
   \en{propo}

   Let us now conclude the proof of \theo{th:main}.
   Thus  after multiplication of \eqre{24071816h02}  by $1+h$, we have  
  \beq  h(1+h)&=&\ell(1+h)+  z^2ab(1-\ell) +  h\ell (1+h)+O\lf(\ff{T\, (\Im z)^2}\ri)
 \eeq
 Using the fact, following from lemma \re{lem:SubGaussianMainVariables}, Cauchy-Schwarz inequality and first part of Proposition \re{247153}, that $$\E \Tta- \f{z^2 ab}{1+h}=O\lf(\ff{T|\Im z|^5}\ri).$$
 We get, for $\tta:=\f{z^2 ab}{1+h}$,
  \beq  h &=&\ell  +  \tta(1-\ell) +  h\ell  +O\lf(\ff{T\, |\Im z|^5}\ri)
 \eeq
 \ie \be\la{130818final}\ell=\f{h-\tta}{1+h-\tta}+O\lf(\ff{T\, |\Im z|^5}\ri).\ee Then, conclude that \eqre{eq:1508188h} is true using \lemm{lem:SubGaussianMainVariables}.

  \subsection{Proof of Theorem \re{th:main2}}
 In the case where $\cA=I_n$ and $\cB=I_p$, the random variables $\wtA$ and $\wtB$ from \eqre{Notation:Intro2507wtAwtB} are respectively equal  to $G$ and $ (p-n)/(Tz^{2})+G$, 
and rather than using \eqre{130818final} to estimate $\ell$, 
we shall solve  \eqre{24071816h02} without using $a$ and $b$.
 Using \eqre{ODODOD} and \eqre{ODODODbis},   after multiplication by $1-\ell$, \eqre{24071816h02} rewrites 
  \beq  h(1-\ell)&=&\ell (1-\ell)+  z^2 g( (p-n)/(Tz^{2})+g)(1+h) +  h\ell (1-\ell)+O\lf(\ff{T\, (\Im z)^2}\ri)
 \eeq
 for $g:=\E G$.
For $\ka:=-  g( (p-n)/T +z^2 g)(1+h) $, we get 
\be\la{3012181}\lf(1+h\ri)\ell^2-\lf(1+2h\ri)\ell+h+\ka+O\lf(\ff{T\, |\Im z|^5}\ri)=0\ee
where we have used the fact, following from lemma \re{lem:SubGaussianMainVariables}, Cauchy-Schwarz inequality and first part of Proposition \re{247153}, that $$\E K- \lf(\f{p-n}T+z^2g\ri)g(1+h)^2=O\lf(\ff{T(\Im z )^5}\ri).$$
Second order polynomial equation \eqre{3012181} solves as  \beq \ell&=&\f{1+2h\pm\sqrt{1-4\ka(1+h)}}{2(1+h)}+O\lf(\ff{T\, |\Im z|^5}\ri)
 \eeq
 Considering the case where $\al$ and $\bet$ are small (where we should have $\ell\approx h$, as explained in \rmq{rmk:boundednp}) and using analytic continuation, we have 
 \be\la{130818final2}\ell=\f{1+2h-\sqrt{1-4\ka(1+h)}}{2(1+h)}+O\lf(\ff{T\, |\Im z|^5}\ri)\ee
 for $\sqrt{\cdote}$ the analytic version of the square root on $\C\bck(-\infty, 0]$ with value $1$ at $1$. Then, again, conclude  using \lemm{lem:SubGaussianMainVariables}.

   \subsection{Proof of Proposition \re{propo:key_relations}}\la{proofpropo:key_relations}
   \subsubsection{Proof of \eqre{24071816h02}: expansion of $\E \Tr \bG\bE\bE'$} 
Using $\bE=\ff T\sum_t X(t)Y(t)'$, we have 
\beqy\nober  \Tr \bE\bE'\bG(z)&=&\ff T\sum_{t=1}^T\Tr X(t)Y(t)'\bE'\bG\\
\nober&=&\ff T\sum_{t=1}^T  Y(t)'\bE'\bG X(t)\\
\nober&=&\ff{2T}\sum_{t=1}^T  Y(t)'\bE'\bG X(t)+X(t)'\bG \bE Y(t)\\
\la{eq:iamilwy}&=&\ff{2T}\sum_{t=1}^T  Z(t)'\bpm 0&\bG\bE\\ \bE'\bG  &0\epm Z(t) 
\eeqy
for \be\la{def:Z(t)}Z(t):=\bpm X(t)\\ Y(t)\epm.\ee

By \eqre{2307181} from \co{cor91171}, it follows that    \begin{align}&
\E\Tr \bE\bE'\bG (z)= \la{eq:20071818h}\\ &\ff 2 \E \Tr \Si \bpm 0&\bG \bE\\ \bE'\bG  &0\epm+\ff{2T}\E\sum_{t=1}^T \sum_{k=1}^{ m}\lf(\Si  \lf( \f{\pa}{\pa Z(t)_k}  \bpm 0&\bG \bE\\ \bE'\bG  &0\epm\ri) Z(t) \ri)_{k }\nober
\end{align}
To distinguish between the $X$ components of $Z(t)$ (the $n$ first ones) and the $Y$ components (the $p$ last ones), we shall now rewrite the above sum as follows: for any $t$, \beqy\la{20071816h41} &&\sum_{k=1}^{ m}\lf(\Si  \lf( \f{\pa}{\pa Z(t)_k}  \bpm 0&\bG \bE\\ \bE'\bG  &0\epm\ri) Z(t)\ri)_{k }\\ \nober
&=&\sum_{k=1}^{ n}\bee_k'\Si  \lf( \f{\pa}{\pa X(t)_k}  \bpm 0&\bG \bE\\ \bE'\bG  &0\epm\ri) Z(t) + 
\sum_{k=1}^{ p}\bee_{n+k}'\Si \lf( \f{\pa}{\pa Y(t)_k}  \bpm 0&\bG \bE\\ \bE'\bG  &0\epm\ri) Z(t) ,
\eeqy
where the $\bee_i$'s denote the (column) vectors of the canonical basis of $\R^m$.

Let us now introduce  the $m\ti m$ matrix \be\la{eq:defF}\bF:=\bpm 0& \bE\\ \bE'&0\epm.\ee
Note that for any $k\ge 0$ integer,  for we have $$\bF^{2k}=\bpm (\bE\bE')^k&0\\ 0&(\bE'\bE)^k\epm,\qquad \bF^{2k+1}=\bpm0& (\bE\bE')^k\bE\\ (\bE'\bE)^k\bE'&0\epm$$
so that for $|z|$ large enough, \beqy\nober (z-\bF)^{-1}&=&\sum_{k\ge 0}\f{\bF^k}{z^{k+1}}\\ \nober&=&\sum_{k\ge 0}z^{-(2k+1)}\bpm (\bE\bE')^k&0\\ 0&(\bE'\bE)^k\epm \\ \nober&&+\sum_{k\ge 0}z^{-2(k+1)}\bpm0& (\bE\bE')^k\bE\\ \bE'(\bE\bE')^k&0\epm\\
&=&\bpm z(z^2-\bE\bE')^{-1}&(z^2-\bE\bE')^{-1}\bE
\\
\bE'(z^2-\bE\bE')^{-1}&z(z^2-\bE'\bE)^{-1}
\epm\nober\\
&=&\bpm z\bG &\bG \bE
\\
\bE'\bG &z\wtG
\epm,  
\la{19071812}
\eeqy 
which is true for all $z\in\C\bck\R$, by analytic continuation. 

\blem For any $t$, we have, for $k=1, \ld, n$,  \beq  \f{\pa}{\pa X(t)_k}  \bpm 0&\bG \bE\\ \bE'\bG  &0\epm&=
&\ff{2T}(z-\bF)^{-1}\bpm 0& \bee_kY(t)'\\
 Y(t)\bee_k'&0\epm (z-\bF)^{-1}+\\ &&
\ff{2T}(z+\bF)^{-1}\bpm 0&\bee_kY(t)'\\ &\\
 Y(t)\bee_k'&0\epm (z+\bF)^{-1}
\eeq
for $\bee_k$ the $k$-th (column) vector of the canonical basis in $\R^n$ and we have, for $k=1, \ld, p$,  \beq \f{\pa}{\pa Y(t)_k}  \bpm 0&\bG \bE\\ \bE'\bG  &0\epm&=&
\ff{2T}(z-\bF)^{-1}\bpm 0&X(t)\bee_k' \\
\bee_kX(t)'&0\epm (z-\bF)^{-1}+\\
&&\ff{2T}(z+\bF)^{-1}\bpm 0& X(t)\bee_k'\\ &\\
 \bee_kX(t)'&0\epm (z+\bF)^{-1}
\eeq
for $\bee_k$ the $k$-th (column) vector of the canonical basis in $\R^p$ 
\elem
\bpr We define  the function $$\vfi(s):=\f{s}{z^2-s^2}=\ff2\lf(\ff{z-s}-\ff{z+s}\ri)\qquad (s\in \R).$$  It is easy to see, by \eqre{19071812},  that we have $$ \bpm 0&\bG \bE\\ \bE'\bG  &0\epm=\vfi(\bF)$$

We want to compute the derivatives,  at $Z(1), \ld,\widehat{Z(t)},\ld,  Z(T)$ fixed, of the function 
$$Z(t)\mapsto \vfi(\bF).$$
The differential of the function $M\mapsto (z-M)^{-1}$ at the matrix $M$ is the operator $$H\mapsto 
 (z-M)^{-1}H (z-M)^{-1},$$the differential of the function $M\mapsto (z+M)^{-1}$ at the matrix $M$ is the operator $$H\mapsto 
- (z+M)^{-1}H (z+M)^{-1},$$
hence the differential of the function $M\mapsto \vfi(M)$ at the matrix $M$ is the operator $$H\mapsto 
\ff2\lf( (z-M)^{-1}H (z-M)^{-1}+(z+M)^{-1}H (z+M)^{-1}\ri).$$ Besides, at $Z(1), \ld,\widehat{Z(t)},\ld,  Z(T)$ fixed, the differential of the function $Z(t)\mapsto \bE$ at  $Z(t)$ is the operator $$\bpm x\\ y\epm\mapsto 
\ff T\lf(X(t)y'+xY(t)'\ri),$$ so that  the differential of the function $Z(t)\mapsto \bF$ at  $Z(t)$ is the operator $$\bpm x\\ y\epm\mapsto 
\ff T\bpm 0&X(t)y'+xY(t)'\\
yX(t)'+Y(t)x'&0\epm.$$ 
It follows that at $Z(1), \ld,\widehat{Z(t)},\ld,  Z(T)$ fixed, the differential of the function $Z(t)\mapsto \vfi(\bF)$ at  $Z(t)$ is the operator  \begin{align*} \bpm x\\ y\epm\mapsto& 
\ff{2T}(z-\bF)^{-1}\bpm 0&X(t)y'+xY(t)'\\
yX(t)'+Y(t)x'&0\epm (z-\bF)^{-1}+\\
&\ff{2T}(z+\bF)^{-1}\bpm 0&X(t)y'+xY(t)'\\ &\\
yX(t)'+Y(t)x'&0\epm (z+\bF)^{-1}\end{align*} The conclusion follows.
 \epr

We deduce that 
  \beq  \f{\pa}{\pa X(t)_k}  \bpm 0&\bG \bE\\ \bE'\bG  &0\epm Z(t)&=
&\ff{2T}(z-\bF)^{-1}\bpm 0& \bee_kY(t)'\\
 Y(t)\bee_k'&0\epm (z-\bF)^{-1}Z(t)+\\ &&
\ff{2T}(z+\bF)^{-1}\bpm 0&\bee_kY(t)'\\ &\\
 Y(t)\bee_k'&0\epm (z+\bF)^{-1}Z(t)
\eeq
By \eqre{19071812},  
\begin{align*}& (z-\bF)^{-1}\bpm 0& \bee_kY(t)'\\
 Y(t)\bee_k'&0\epm (z-\bF)^{-1}=\\
 &
 \bpm z\bG &\bG \bE
\\
\bE'\bG &z\wtG
\epm
\bpm 0& \bee_kY(t)'\\
 Y(t)\bee_k'&0\epm
 \bpm z\bG &\bG \bE
\\
\bE'\bG &z\wtG
\epm=\\
&\bpm
z\bG  \bee_kY(t)' \bE'\bG   + z\bG \bE Y(t)\bee_k' \bG   & 
z^2\bG  \bee_kY(t)' \wtG  + \bG \bE Y(t)\bee_k' \bG \bE  \\
\bE'\bG  \bee_kY(t)' \bE'\bG   + z^2\wtG Y(t)\bee_k' \bG   & 
z\bE'\bG  \bee_kY(t)' \wtG  + z\wtG Y(t)\bee_k' \bG \bE 
\epm 
 \end{align*}
Thus 
\begin{align*}& (z-\bF)^{-1}\bpm 0& \bee_kY(t)'\\
 Y(t)\bee_k'&0\epm (z-\bF)^{-1}Z(t)=\\
&\bpm
z\bG  \bee_kY(t)' \bE'\bG  X(t) + z\bG \bE Y(t)\bee_k' \bG   X(t)+
z^2\bG  \bee_kY(t)' \wtG Y(t)  + \bG \bE Y(t)\bee_k' \bG \bE  Y(t)\\
\bE'\bG  \bee_kY(t)' \bE'\bG  X(t) + z^2\wtG Y(t)\bee_k' \bG X(t)  +
z\bE'\bG  \bee_kY(t)' \wtG  Y(t)+ z\wtG Y(t)\bee_k' \bG \bE Y(t)
\epm 
 \end{align*}
  Then, it is easy to see, by \eqre{19071812}, that computing $$(z+\bF)^{-1}\bpm 0& \bee_kY(t)'\\
 Y(t)\bee_k'&0\epm (z+\bF)^{-1}Z(t)$$ amounts to take the same formula and change $\bE$ into $-\bE$.  
 After, adding both and dividing by $2T$ amounts to    keep only,  in the previous formula,  the terms which are even in $\bE$ (and divide by $T$). We get
\be\la{derX2007}  \f{\pa}{\pa X(t)_k}  \bpm 0&\bG \bE\\ \bE'\bG  &0\epm Z(t)= 
\ff{T}\bpm z^2\bG  \bee_kY(t)' \wtG Y(t) + \bG \bE  Y(t)\bee_k' \bG \bE Y(t)\\ 
 z^2\wtG Y(t)\bee_k' \bG X(t)+ \bE'\bG  \bee_kY(t)' \bE'\bG X(t)  \epm
 \ee

In the same way, 
\be\la{derY2007}  \f{\pa}{\pa Y(t)_k}  \bpm 0&\bG \bE\\ \bE'\bG  &0\epm Z(t)= 
\ff{T}\bpm z^2\bG  X(t)\bee_k' \wtG Y(t) + \bG \bE  \bee_kX(t)' \bG \bE Y(t)\\ 
 z^2\wtG \bee_kX(t)' \bG X(t)+ \bE'\bG  X(t)\bee_k' \bE'\bG X(t)  \epm
 \ee

 Let us write $$\Si=\bpm\cA&\cC\\ \cC'&\cB\epm, \quad \cA=\Cov(X),\; \cC=\Cov(X,Y), \; \cB=\Cov(Y).$$
 
 By \eqre{20071816h41} , \eqre{derX2007} and \eqre{derY2007} (and using \eqre{eq:Trace1}, \eqre{eq:Trace2} and the facts that $\bG '=\bG $ and $\wtG'=\wtG$), we have 
  \beqy\nober &&\sum_{k=1}^{ m}\lf(\Si  \lf( \f{\pa}{\pa Z(t)_k}  \bpm 0&\bG \bE\\ \bE'\bG  &0\epm\ri) Z(t)\ri)_{k }\\ 
  \nober&=&\ff T\sum_{k=1}^n \bee_{k }'\cA\lf(z^2\bG  \bee_{k }Y(t)' \wtG Y(t) + \bG \bE  Y(t)\bee_{k }' \bG \bE Y(t) \ri)+\\
   \nober&&\ff T\sum_{k=1}^n\bee_{k }'\cC\lf( z^2\wtG Y(t)\bee_{k }' \bG X(t)+ \bE'\bG  \bee_{k }Y(t)' \bE'\bG X(t) \ri)+\\ 
   \nober&& \ff T\sum_{k=1}^p \bee_{k }'\cC'\lf( z^2\bG  X(t)\bee_{k }' \wtG   Y(t) + \bG \bE  \bee_{k }X(t)' \bG \bE Y(t) \ri)+\\
    \nober&&\ff T\sum_{k=1}^p\bee_{k }'\cB\lf( z^2\wtG   \bee_{k }X(t)' \bG X(t)+ \bE'\bG  X(t)\bee_{k }' \bE'\bG X(t)\ri)\\ 
  \nober&=&\ff T\lf(z^2(\Tr \cA \bG ) Y(t)'\wtG   Y(t)+Y(t)'\bE'\bG \cA \bG \bE  Y(t)\ri)+\\
\nober&&\f2T\lf(z^2X(t)'\bG \cC \wtG   Y(t) + (\Tr \cC\bE'\bG  ) Y(t)' \bE'\bG X(t) \ri)+\\ 
 \la{20071816h4177}&&\ff T\lf(z^2(\Tr \cB \wtG  ) X(t)'\bG  X(t)+X(t)'\bG \bE \cB \bE'\bG   X(t)\ri)
    \eeqy
    Let us now sum \eqre{20071816h4177} over $t=1, \ld, T$. Having in mind that      $$ \ff T \sum_{t=1}^T X(t)Y(t)'=\bE,\qquad \bC_X:= \ff T \sum_{t=1}^T X(t)X(t)'\AND 
    \bC_Y:= \ff T \sum_{t=1}^T Y(t)Y(t)',$$ we get 
     \begin{align} &  \sum_{t=1}^T\sum_{k=1}^{ m}\lf(\Si  \lf( \f{\pa}{\pa Z(t)_k}  \bpm 0&\bG \bE\\ \bE'\bG  &0\epm\ri) Z(t)\ri)_{k }=\nober\\ \nober&z^2 \Tr \cA \bG   \Tr  \wtG   \bC_Y+\Tr \bE'\bG \cA \bG \bE  \bC_Y+2z^2\Tr \bG \cC \wtG   \bE' +\\ & 2 \Tr \cC\bE'\bG    \Tr \bE'\bG \bE+
 \la{20071816h4177sumed}  z^2 \Tr \cB \wtG   \Tr \bG  \bC_X+\Tr \bG \bE \cB \bE'\bG \bC_X    \end{align}
    Joining \eqre{eq:20071818h} and \eqre{20071816h4177sumed}, we get 
   \begin{align}&
\E\Tr \bE\bE'\bG (z)=\nober\\
\nober&  \E \Tr  \cC\bE'\bG  +\ff{2T}\E \Big[z^2 \Tr \cA \bG   \Tr  \wtG   \bC_Y+\Tr \bE'\bG \cA \bG \bE  \bC_Y+2z^2\Tr \bG \cC \wtG   \bE' + \nober\\ &2 \Tr \cC\bE'\bG    \Tr \bE'\bG \bE+
 \la{StMaur210718}  z^2 \Tr \cB \wtG   \Tr \bG  \bC_X+\Tr \bG  \bE\cB \bE'\bG \bC_X   \Big],
\end{align}
which allows to conclude.

   \subsubsection{Proof of \eqre{ODODOD}: expansion of $ \E  \Tr \bG \bC_X$}
For $\bF$ as in \eqre{eq:defF}, by \eqre{19071812}, we have \beqy\la{1308181} \bpm  \bG &0\\ 0&  \wtG  \epm &=&\psi(\bF)\eeqy for  \beqy\la{1308182}\psi(s)&:=&\ff{z^2-s^2}\;=\;\ff{2z}\lf(\ff{z-s}+\ff{z+s}\ri).\eeqy
It follows that for $P:=\bpm I_n&0\\ 0&0\epm$, $$\bpm \bG &0\\ 0&0\epm=P\psi(\bF).$$

For $Z(t)$ as defined in \eqre{def:Z(t)}, we have, by \eqre{2307181} of \co{cor91171},  
\beq \E\Tr \bG \bC_X&=&\ff T \sum_{t=1}^T \E X(t)'\bG X(t)\\
&=&\ff T \sum_{t=1}^T \E Z(t)'\bpm \bG &0\\ 0&0\epm Z(t)\\
&=&\E \Tr \cA \bG +\ff T\E \sum_{t=1}^T \sum_{k=1}^m\bee_k'  \Si \f{\pa}{\pa Z(t)_k} \lf(P\psi(\bF)\ri)Z(t)\\
&=&\E \Tr \cA \bG +\\
&&\ff T\E \sum_{t=1}^T\lf( \sum_{k=1}^n\bee_k'  \Si \f{\pa}{\pa X(t)_k} \lf(P\psi(\bF)\ri)Z(t)+ \sum_{l=1}^p\bee_{n+l}'  \Si \f{\pa}{\pa Y(t)_l} \lf(P\psi(\bF)\ri)Z(t)\ri) \eeq
Note that \beq \f{\pa}{\pa X(t)_k} P\psi(\bF) &=&\ff{2zT}P(z-\bF)^{-1}\bpm 0&\bee_kY(t)'\\ Y(t)\bee_k'&0\epm(z-\bF)^{-1}\\ &&-\ff{2zT}P(z+\bF)^{-1}\bpm 0&\bee_kY(t)'\\ Y(t)\bee_k'&0\epm(z+\bF)^{-1}\\
&=&\ff{T}\bpm
\bG  \bee_kY(t)' \bE'\bG   + \bG \bE Y(t)\bee_k' \bG &0\\ 0&0\epm
\eeq
and in the same way, 
\beq \f{\pa}{\pa Y(t)_l} P\psi(\bF)  
&=&\ff{T}\bpm
\bG  X(t)\bee_l' \bE'\bG   + \bG \bE  \bee_lX(t)' \bG &0\\ 0&0\epm
\eeq
We deduce that 
\beq \E\Tr \bG \bC_X
&=&\E \Tr \cA \bG +\ff{T^2}\E \sum_{t=1}^T  \sum_{k=1}^n\bee_k'  \cA\lf(\bG  \bee_kY(t)' \bE'\bG   + \bG \bE Y(t)\bee_k' \bG \ri)X(t)\\
&&+\ff{T^2}\E \sum_{t=1}^T  \sum_{l=1}^p\bee_l'  \cC'\lf( \bG  X(t)\bee_l' \bE'\bG   + \bG \bE  \bee_lX(t)' \bG \ri)X(t)\\
&=&\E \Big[\Tr \cA \bG +\ff T
\Tr \cA \bG \Tr \bG \bE\bE'+\ff{T^2} \sum_{t=1}^T   X(t)'\bG \cA \bG \bE Y(t)\\
&&+\ff{T^2}  \sum_{t=1}^T     X(t)'\bG \bE \cC'\bG  X(t)   + \ff T\Tr \cC'\bG \bE   \Tr \bG \bC_X\Big]\\
&=&\E \Big[\Tr \cA \bG +\ff T
\Tr \cA \bG \Tr \bG \bE\bE'+\ff{T}\Tr \bG \cA \bG \bE\bE'\\
&&+\ff{T} \Tr \bG \bE \cC'\bG  \bC_X   + \ff T\Tr \cC'\bG \bE   \Tr \bG \bC_X\Big],
\eeq
which allows to conclude.

      \subsubsection{Proof of \eqre{ODODODbis}: expansion of $ \E  \Tr \wtG   \bC_Y$}
For $\bF$ as in \eqre{eq:defF}, by \eqre{19071812}, we have \beq \bpm  \bG &0\\ 0&  \wtG  \epm &=&\psi(\bF)\eeq for $$\psi(s):=\ff{z^2-s^2}=\ff{2z}\lf(\ff{z-s}+\ff{z+s}\ri).$$
It follows that for $Q:=\bpm 0&0\\ 0&I_p\epm$, $$\bpm 0&0\\ 0&\wtG  \epm=P\psi(\bF).$$

For $Z(t)$ as defined in \eqre{def:Z(t)}, we have, by \eqre{2307181}  of \co{cor91171},  
\beq \E\Tr \wtG  \bC_Y&=&\ff T \sum_{t=1}^T \E Y(t)'\wtG   Y(t)\\
&=&\ff T \sum_{t=1}^T \E Z(t)'\bpm 0&0\\ 0&\wtG \epm Z(t)\\
&=&\E \Tr \cB \wtG  +\ff T\E \sum_{t=1}^T \sum_{k=1}^m\bee_k'  \Si \f{\pa}{\pa Z(t)_k} \lf(P\psi(\bF)\ri)Z(t)\\
&=&\E \Tr \cB \wtG  +\ff T\E \sum_{t=1}^T\Bigg( \sum_{k=1}^n\bee_k'  \Si \f{\pa}{\pa X(t)_k} \lf(P\psi(\bF)\ri)Z(t)+ \\&&\qquad\qquad\qquad\qquad\qquad\qquad \sum_{l=1}^p\bee_{n+l}'  \Si \f{\pa}{\pa Y(t)_l} \lf(P\psi(\bF)\ri)Z(t)\Bigg) \eeq
Note that \beq \f{\pa}{\pa X(t)_k} P\psi(\bF) &=&\ff{2zT}P(z-\bF)^{-1}\bpm 0&\bee_kY(t)'\\ Y(t)\bee_k'&0\epm(z-\bF)^{-1}\\ &&-\ff{2zT}P(z+\bF)^{-1}\bpm 0&\bee_kY(t)'\\ Y(t)\bee_k'&0\epm(z+\bF)^{-1}\\
&=&\ff{T}\bpm
0&0\\ 0&\bE'\bG  \bee_kY(t)' \wtG    + \wtG   Y(t)\bee_k' \bG \bE \epm
\eeq
and in the same way, 
\beq \f{\pa}{\pa Y(t)_l} P\psi(\bF)  
&=&\ff{T}\bpm
0&0\\ 0&\bE'\bG  X(t)\bee_l' \wtG    + \wtG   \bee_lX(t)' \bG \bE \epm\eeq
We deduce that 
\beq  \E\Tr \wtG  \bC_Y
&=&\E \Tr \cB \wtG  +
\ff{T^2}\E \sum_{t=1}^T  \sum_{k=1}^n\bee_k'  \cC\lf(\bE'\bG  \bee_kY(t)' \wtG    + \wtG   Y(t)\bee_k' \bG \bE\ri)Y(t)\\
&&+\ff{T^2}\E \sum_{t=1}^T  \sum_{l=1}^p\bee_l'  \cB\lf( \bE'\bG  X(t)\bee_l' \wtG    + \wtG   \bee_lX(t)' \bG \bE\ri)Y(t)\\
&=&\E \Big[\Tr \cB \wtG  +\ff T
\Tr \cC\bE'\bG  \Tr \wtG  \bC_Y+ \ff{T^2} \sum_{t=1}^T Y(t)'\bE'\bG \cC  \wtG   Y(t)\\
&&+\ff{T^2} \sum_{t=1}^T Y(t)'\wtG  \cB \bE'\bG  X(t)+\ff T\Tr \cB \wtG  \Tr \bG \bE\bE'
\Big]\\
&=&\E \Big[\Tr \cB \wtG  +\ff T
\Tr \cC\bE'\bG  \Tr \wtG  \bC_Y+ \ff{T}\Tr \bE'\bG \cC  \wtG   \bC_Y\\
&&+\ff{T}\Tr \wtG  \cB \bE'\bG  \bE+\ff T\Tr \cB \wtG  \Tr \bG \bE\bE'
\Big],
\eeq which allows to conclude.

      \subsection{Proof of Proposition \re{prop_unbiased_estimator}}

  \subsubsection{Proof of \eqre{precise_estimate_cleaned_true_emp1}}
  A simple application of equality $s\ucleaned_k=\bu_k'\cC\bv_k$ from \eqre{0507181} gives: 
  $$\E  \sum_{k=1}^n s_ks_k\ucleaned=\E \sum_{k=1}^n s_k\bu_k'\cC\bv_k=\E \Tr \bE'\cC= \Tr (\E\bE)'\cC= \Tr \cC'\cC= \sum_{k=1}^n (s_k\utrue)^2.$$
    \subsubsection{Proof of \eqre{precise_estimate_cleaned_true_emp2} and \eqre{precise_estimate_cleaned_true_emp3}} We start with the following Gaussian integrals:
    \blem\la{lem:FreeProbaOrder4} Let $m,T\ge 1$, $P,Q\in\R^{m\ti m}$  and let $Z\in\R^{m\ti T}$ be a matrix whose entries  are independent standard Gaussian variables.
Then   we have \be\la{301218331}\E\Tr ZZ'QZZ'P = T^2\Tr PQ+T\Tr P'Q+T\Tr P\Tr Q \ee
and  \be\la{30121833}\E\Tr ZZ'Q\Tr ZZ'P = T^2\Tr P\Tr Q +T\Tr PQ+T\Tr P'Q.\ee
\elem 

\bpr We have 
\beq\E\Tr ZZ'QZZ'P
&=&\E\sum_{i,j,k,l,r,s }\E Z_{ij}Z_{kj}Q_{kl}Z_{lr}Z_{sr}P_{si}
\eeq
Using then  the fact that the entries of $Z$ are even and independent, we get 
\beq\E \Tr ZZ'QZZ'P&=&
\sum_{i,j,l,r}\E Z_{ij}Z_{ij}Q_{il}Z_{lr}Z_{lr}P_{li}
+
\sum_{i,j,k} \E Z_{ij}Z_{kj}Q_{ki}Z_{ij}Z_{kj}P_{ki}
+\\ &&
\sum_{i,j,k} \E Z_{ij}Z_{kj}Q_{kk}Z_{kj}Z_{ij}P_{ii}
-2\sum_{i,j}\E Z_{ij}Z_{ij}Q_{ii}Z_{ij}Z_{ij}P_{ii}\\
&=&
\sum_{(i,j)\ne (l,r)}Q_{il}P_{li}
+
3\sum_{i,j}  Q_{ii}P_{ii}
+ 
\sum_{i\ne k,j}  Q_{ki} P_{ki}
+
3\sum_{i,j} Q_{ii}P_{ii}
+\\ &&
\sum_{i\ne k,j}  Q_{kk} P_{ii}
+
3\sum_{i,j}  Q_{ii} P_{ii}
-6\sum_{i,j} Q_{ii} P_{ii}\\
&=&
\sum_{i,j,l,r}Q_{il}P_{li}-\sum_{i,j}Q_{ii}P_{ii}
+
3\sum_{i,j}  Q_{ii}P_{ii}
+\\ &&
\sum_{i, k,j}  Q_{ki} P_{ki}
-\sum_{i,j}  Q_{ii} P_{ii}
+
3\sum_{i,j} Q_{ii}P_{ii}
+\\ &&
\sum_{i, k,j}  Q_{kk} P_{ii}
-\sum_{i ,j}  Q_{ii} P_{ii}
+
3\sum_{i,j}  Q_{ii} P_{ii}
-6\sum_{i,j} Q_{ii} P_{ii}\\
&=&T^2\Tr PQ+T\Tr P'Q+T\Tr P\Tr Q 
\eeq
In the same way,
 \beq\E\Tr ZZ'Q\Tr ZZ'P
&=&\E\sum_{i,j,k,l,r,s }\E Z_{ij}Z_{kj}Q_{ki}Z_{lr}Z_{sr}P_{sl}
\eeq
and using again  the fact that the entries of $Z$ are even and independent, we get 
\beq\E \Tr ZZ'Q\Tr ZZ'P&=&
\sum_{i,j,l,r}\E Z_{ij}Z_{ij}Q_{ii}Z_{lr}Z_{lr}P_{ll}
+
\sum_{i,j,k} \E Z_{ij}Z_{kj}Q_{ki}Z_{ij}Z_{kj}P_{ki}\\ &&
+
\sum_{i,j,k} \E Z_{ij}Z_{kj}Q_{ki}Z_{kj}Z_{ij}P_{ik}
-2\sum_{i,j}\E Z_{ij}Z_{ij}Q_{ii}Z_{ij}Z_{ij}P_{ii}\\
&=&
\sum_{(i,j)\ne(l,r)} Q_{ii} P_{ll}+3\sum_{i,j } Q_{ii} P_{ii}
+
\sum_{i\ne k ,j}  Q_{ki} P_{ki}
+
3\sum_{i,j}  Q_{ii} P_{ii}\\ &&
+
\sum_{i\ne k,j}  Q_{ki} P_{ik}
+
3\sum_{i,j } Q_{ii} P_{ii}
-6\sum_{i,j} Q_{ii} P_{ii}
\\
&=&
\sum_{i,j,l,r} Q_{ii} P_{ll}
-\sum_{i,j } Q_{ii} P_{ii}
+3\sum_{i,j } Q_{ii} P_{ii}
+
 \sum_{i, k ,j}  Q_{ki} P_{ki}- \sum_{i ,j}  Q_{ii} P_{ii}
\\ &&
+
\sum_{i,k,j}  Q_{ki} P_{ik}- \sum_{i ,j}  Q_{ii} P_{ii}
\\
&=&
\sum_{i,j,l,r} Q_{ii} P_{ll}
+
 \sum_{i, k ,j}  Q_{ki} P_{ki} 
+
\sum_{i,k,j}  Q_{ki} P_{ik} \\
&=&
  T^2\Tr P\Tr Q +T\Tr PQ+T\Tr P'Q
\eeq
\epr

Let us now prove \eqre{precise_estimate_cleaned_true_emp2} and \eqre{precise_estimate_cleaned_true_emp3}. 
Recall that  $$\bbm X\\ 
Y\ebm\sim\cN(0, \Si)$$ for  $\Si=\bpm\cA & \cC\\ \cC'&\cB\epm\in \R^{m\ti m}$ for $m=n+p$.

Equations \eqre{precise_estimate_cleaned_true_emp2} and \eqre{precise_estimate_cleaned_true_emp3} follow directly from the following lemma.
\blem\la{lem:151121} We have $$\E \Tr \bE\bE'=\f{T+1}T\Tr \cC\cC'+T^{-1}\Tr \cA\Tr \cB$$
and $$ \E \Tr \bC_X\Tr \bC_Y= \Tr \cA\Tr \cB+2T^{-1}\Tr \cC\cC'$$
\elem

\bpr Let  $\bZ\in\R^{m\ti T}$ be a matrix  whose entries  are independent standard Gaussian variables. Let   $R  \in\R^{n\ti m}$ and $S \in\R^{p\ti m}$ so that $\bX$,  $\bY$ can be realized by 
  $$\bX=R\bZ\AND \bY=S\bZ.$$ Then by \eqre{defbE} and \eqre{eq:defCXCY}:
\be\la{eq:30121822}\bC_X=\ff T\bX\bX'=\ff TR\bZ\bZ'S',\; \bC_Y=\ff T\bY\bY'=\ff T S\bZ\bZ'S,\; \bE=\ff T\bX\bY'=\ff TR\bZ\bZ'S'.\ee
 
 By Lemma \re{lem:FreeProbaOrder4}, we have, for $P:=R'R$ and $Q:=S'S$, 
\beq T^2\E \Tr \bE\bE'&=& \E\Tr R\bZ\bZ'S'S\bZ\bZ'R'
\\
&=& \E\Tr \bZ\bZ'Q\bZ\bZ'P 
\eeq
and 
\beq T^2\E \Tr \bC_X\Tr \bC_Y&=& \E\Tr R\bZ\bZ'R'\Tr S\bZ\bZ'S'
\\
&=& \E\Tr \bZ\bZ'Q\Tr \bZ\bZ'P 
\eeq
Then, we conclude  noting that $\cC=RS'$, $\cA=RR'$ and $\cB=SS'$.
\epr

  \subsection{Proof of concentration results}
      \subsubsection{Proof of \lemm{lem:SubGaussianMainVariables}}\la{sec:measure_concentration}
     With the notation of the proof of Lemma  \re{lem:151121}, by 
     the second part of Proposition \re{247153}, what we have to prove is that the functions mapping  $\bZ\in \R^{m\ti T}$  to   the variables $G,H,A,B,L,\Tta,K\in \C$  are all Lipschitz for the Frobenius norm $\|\cdot\|_{\text{F}}$  from \eqre{opt_problem00} on $\R^{m\ti T}$, 
      with Lipschitz constant  $$O\lf(\ff{T(\Im z )^4}\ri).$$As this argument is quite standard (close to e.g. \cite[Sec. 2.3.1]{agz} or \cite[Lem.  7.1]{capJTP} with \cite[Lem. B.2]{FloRomain1} instead of \cite[Lem.  A.2]{capJTP}), we only give the main lines. Consider a variation $\delta_\bZ$ of $\bZ$, and then:
      \begin{enumerate}\ite Use \eqre{eq:30121822}.
            \ite Use the \emph{resolvant formula}: for all square matrices $M,\delta_M$, $$(z-(M+\delta_M))^{-1}-(z-M)^{-1}\;=\;(z-(M+\delta_M))^{-1}\delta_M(z-M)^{-1}$$  to expand the variations of the matrices $\bG$ and $\wtG$ at first order in $\delta_\bZ$.
      \ite By non-commutative H\"older inequalities (see e.g. \cite[Appendix A.3]{agz}), for any product $M_1\cd M_k$ of matrices with any size and any $i=1,\ld, k$, $$\|M_1\cd M_k\|_{\text{F}}\le \|M_1\|_{\op{op}}\cd \widehat{\|M_{i}\|}_{\op{op}}\cd \|M_k\|_{\op{op}} \|M_{i}\|_{\text{F}}$$  where  $\|\cdot\|_{\op{op}}$ denotes the operator norm. This has to be used with the fact that $\bG$ and $\wtG$ have operator norms $\le \op{dist}(z^2, [0, +\infty))^{-1}$.
      \ite On any square matrices space endowed with the Frobenius norm, the trace is the scalar product with the identity matrix, hence is Lipschitz with Lipschitz constant the Frobenius norm of the identity matrix (which depends on the dimension).
      \end{enumerate}
      
      \subsubsection{Concentration lemma for Section \re{sec:overfiiting}}
      
      \blem\la{lemma_concentration_overfitting} With the notation from Section \re{subsec:modelresultsalgos}, for any deterministic vectors $u\in \R^n$, $v\in \R^p$, the random variable $u'\bE v-u'\cC v$ is centered with $L^2$-norm $\le \sqrt2\|\Si\|_{\op{op}}/\sqrt{T}$.      \elem
      
      \brem Using Hanson-Wright inequality \cite{Vershynin}, one could improve the variance control up  to an exponential control on the tail.\erem
    
      \bpr With the notation of the proof of Lemma  \re{lem:151121},  \be\la{30121833222}\E u'\bE v=u'RS' v=u'\cC v.
      \ee
   Secondly, we have   
  $$u'\bE v=   \ff Tu'R\bZ\bZ'S' v=\ff T\Tr \bZ\bZ'S' v u'R$$
  so that, by \eqre{30121833}, 
\beq \E (u'\bE v)^2&=& (\Tr S' v u'R)^2+\ff T\Tr S' v u'RS' v u'R+\ff T\Tr R'u  v'SS' v u'R\\
&=&
 (u'RS' v )^2+\ff T(u'RS' v)^2+\ff T(u'R R'u)(  v'SS' v ),
\eeq
which, by \eqre{30121833222}, allows to conclude.
      \epr
      
       \section{Appendix}\subsection{Stieltjes transform inversion}Any 
    signed measure $\mu$ on $\R$ can be recovered out of its Stieltjes transform\be\la{def:Stieltjes}g_\mu(z):=\int\f{\ud\mu(t)}{z-t} , \quad z\in \C\bck\R\ee by the formula   \be\la{muoutofg}\mu=-\ff\pi\lim_{\eta\to 0^+}(\Im g_\mu(x+\ii\eta)\ud x),\ee where the limit holds in the weak topology   (see e.g. \cite[Th. 2.4.3]{agz} and use the decomposition of any signed measure as a difference of finite  positive measures). 
    
    Note that for any  $z \in \C\bck\R$ and any $s\ge 0$, we have 
   $$   \ff2\lf(\ff{z-s}-\ff{z+s}\ri)=\f{s}{z^2-s^2}, \qquad \ff2\lf(\ff{z-s}+\ff{z+s}\ri)=\f{z}{z^2-s^2},$$
so that for $s_1,\ld, s_n\ge 0$ and $\rho_1, \ld,\rho_n\in \R$, the  Stieltjes transforms of the measures  $m:=\ff{2n}\sum_{k=1}^n\rho_k(\del_{s_k}-\del_{s_k})$ and $\nu:=\ff{2n}\sum_{k=1}^n\rho_k(\del_{s_k}+\del_{s_k})$ rewrite \be\la{eq:24aout20211}g_m(z)=\ff{n}\sum_{k=1}^n\rho_k \f{s_k}{z^2-s_k^2}\ee and 
 \be\la{eq:24aout20212}g_\nu(z)=\ff{n}\sum_{k=1}^n\rho_k \f{z}{z^2-s_k^2}.\ee

 \subsection{Linear algebra}We notify  some formulas frequently used (and referred to) here:
  for $(\bee_k)$ a collection of column vectors defining an  orthonormal basis,  for any matrices $M,N$, 
\be\la{eq:Trace1}\bee_k' M\bee_l=\bee_l' M'\bee_k,\;\;\;\sum_k \bee_k' M\bee_k=\Tr M,\; \;\;\sum_{k,l}\bee_k' M\bee_l \bee_k' N\bee_l=\Tr MN' \ee and for any column vectors $u,v$, \be\la{eq:Trace2}\sum_k\bee_k' u\bee_k'v=\sum_ku'\bee_k v'\bee_k=v'u.\ee

 \subsection{Stein formula for     Gaussian random vectors}
 \beg{propo}\la{SteinMultidim}Let $X=(X_1, \ld, X_d)$ be a   centered Gaussian vector with covariance $\Si$ and $f : \R^d  \to \R$ be a $\cC^1$    function with derivatives having at most polynomial growth. Then for all $i_0=1, \ld, d$, $$\E X_{i_0}f(X_1, \ld, X_d) \;=\;\sum_{k=1}^d\Si_{i_0k}  \E (\partial_kf)(X_1, \ld, X_d). $$  
 \en{propo} (see e.g. \cite[Lem. A.1]{FloRomain1})

 \bcor\la{cor91171} With the same notation, considering $X$ as a column vector,  for $F:\R^d\to\R^{d\ti d}$ a matrix-valued function, 
 we have   \beqy\la{2307181} \E X'F(X)X &=&\Tr \Si\E F(X)+\sum_{k=1}^d\lf(\E \Si (\pa_kF)(X) X\ri)_k .  
\eeqy  \ecor
 
 \bpr 
 We have, by  Proposition \re{SteinMultidim}, \beq \E X'F(X)X &=&\sum_{ij}\E X_{i}X_{j}F(X)_{ij}\\
  &=&\sum_{ijk}\E \Si_{ik}\f{\pa}{\pa X_k}X_{j}F(X)_{ij}\\
    &=&\sum_{ijk}\E \Si_{ik} \lf(\del_{j=k}F(X)_{ij}+X_{j}(\pa_kF)(X)_{ij}\ri)\\
     &=&\Tr \Si\E F(X)+\sum_{ijk}\E \Si_{ik}  X_{j}(\pa_kF)(X)_{ij} \\
          &=&\Tr \Si\E F(X)+\sum_{k}\lf(\E \Si (\pa_kF)(X) X\ri)_k
 \eeq
 \epr
 
 \subsection{Concentration of measure for Gaussian vectors}The following proposition can be found   e.g. in   \cite[Sec. 4.4.1]{agz} or \cite[Th. 5.2.2]{Vershynin}.
 \beg{propo}\la{247153}Let $X=(X_1, \ld, X_d)$ be a standard real Gaussian vector and $f : \R^d  \to \R$ be a $\mc{C}^1$  function with  gradient $\nabla f$. Then we have \be\la{247151}\op{Var}(f(X))\;\le\; \E \|\nabla f(X)\|^2,\ee where $\|\,\cdot\,\|$ denotes the standard Euclidean norm. 

Besides, if $f$ is $k$-Lipschitz, then for any $t>0$, we have   \be\la{247152}\p (|f(X)-\E f(X)|\ge t)\; \le \; 2\me^{-\f{t^2}{2k^2}},
\ee \ie  $f(X)-\E f(X)$ is Sub-Gaussian with Sub-Gaussian norm $\le k$, up to a universal constant factor. 
 \en{propo}

  \begin{thebibliography}{10}
 \bibitem{popescu}   Amsalu, S.,   Duan, J., Matzinger, H.,  Popescu, I. \emph{Recovery of spectrum from estimated covariance matrices and  statistical kernels for machine learning and big data}, arXiv.

  \bibitem{agz} Anderson, G., Guionnet, A., Zeitouni, O.  \emph{An Introduction to Random Matrices}. Cambridge Studies in Advanced Mathematics, {118} (2009).
  

  \bibitem{FloRIE} Benaych-Georges, F.  \emph{A very short proof of Ledoit-P\'ech\'e's RIE formula for covariance matrices}. Unpublished note available at \url{http://www.cmapx.polytechnique.fr/~benaych/Short_proof_of_Ledoit_Peche.pdf}

  \bibitem{FloRomain1} Benaych-Georges, F., Couillet, R. \emph{Spectral analysis of the Gram matrix of mixture models}, ESAIM Probab. Statist., Vol. 20 (2016), 217--237. 
  
   \bibitem{FloAntti}  Benaych-Georges, F., Knowles, A. \emph{Local semicircle law for Wigner matrices. Advanced topics in random matrices}, 1--90, Panor. Synth\`eses, 53, Soc. Math. France, Paris, 2017.

  \bibitem{bose} Bose, A., Bhattacharjee, M. \emph{Large Covariance and Autocovariance Matrices}, 2018, Chapman and Hall/CRC.

  \bibitem{CorrBouchaud} Bouchaud, J.-P.,  Laloux, L., Miceli, M.A.,  Potters M. \emph{Large dimension forecasting models and random singular value spectra}, Eur. Phys. J. B (2007) 55: 201.

 \bibitem{BBPoRisk} Bun, J., Bouchaud,  J.-P., Potters, M. \emph{Cleaning correlation matrices}, Risk magazine,
2016.
  
 \bibitem{BBPo} Bun, J., Bouchaud,  J.-P., Potters, M. \emph{Cleaning large correlation matrices: Tools from Random Matrix Theory}, Physics Reports
Volume 666, Review article,   1--109, 2017.

  \bibitem{capJTP}  Capitaine, M.  \emph{Additive/multiplicative free subordination property and limiting eigenvectors of spiked additive deformations of Wigner matrices and spiked sample covariance matrices}. J. Theoret. Probab. 26 (2013), no. 3, 595--648.

\bibitem{Robust1} Couillet, R., Pascal, F., Silverstein, J., \emph{Robust estimates of covariance matrices in the large dimensional regime}. IEEE Trans. Inform. Theory 60 (2014), no. 11, 7269--7278.

\bibitem{Robust2} Couillet, R., Pascal, F., Silverstein, J., \emph{The random matrix regime of Maronna's M-estimator with elliptically distributed samples}. J. Multivariate Anal. 139 (2015), 56--78.

\bibitem{ElKaroui}  El Karoui, N. \emph{Spectrum estimation for large dimensional covariance matrices using
random matrix theory}. Annals of Statistics, 2008, 36(6):2757--2790.
\bibitem{ElKaroui2}  El Karoui, N. \emph{Operator norm consistent estimation of large dimensional sparse covariance matrices}, Annals of Statistics, 2008, 36(6): 2717--2756

\bibitem{ElKaroui3}  El Karoui, N. \emph{Random matrices and high-dimensional statistics: beyond covariance matrices}. Proceedings of the International Congress of Mathematicians, 2018.

\bibitem{ErYauS} Erd\H{o}s, L., Schlein, B., Yau, H.-T. \emph{Semicircle law on short scales and delocalization of eigenvectors for Wigner random matrices}. Ann. Probab. 37 (2009), no. 3, 815--852.

\bibitem{ErYau} Erd\H{o}s, L., Yau, H.-T. \emph{A dynamical approach to random matrix theory}. Courant Lecture Notes in Mathematics, 28, New York; American Mathematical Society, Providence, RI, 2017. 

 \bibitem{Tibshi}  Tibshirani, R., Hastie, T. Friedman, J. \emph{Elements of Statistical Learning} Second Edition,  Print 10. Springer Series in Statistics, 2008.

\bibitem{klopp} Klopp, O., Lounici, K., Tsybakov, A. B. \emph{Robust matrix completion}. Probab. Theory Related Fields 169 (2017), no. 1--2, 523--564.

\bibitem{klopp2}   Klopp, O., Y. Lu, A. Tsybakov, A. B., Zhou, H. \emph{Structured Matrix Estimation and Completion}, Bernoulli 25 (4B), 2019, 3883--3911 (2019).

\bibitem{KLT1} Koltchinskii, V., Lounici, K., Tsybakov, A. B. \emph{Estimation of low-rank covariance function}, Stochastic Process. Appl. 126 (2016), no. 12, 3952--3967.

\bibitem{LCBP}  Laloux, L.,   Cizeau, P.,  Bouchaud, J.-P., Potters, M. \emph{Noise Dressing of Financial Correlation Matrices}, 
Phys. Rev. Lett. 83, 1467, 1999.  
\bibitem{LedoitWolf2004}  Ledoit, O., Wolf, M. \emph{A well-conditioned estimator for large-dimensional covariance matrices},  Journal of Multivariate Analysis 88, 2004, 365--411. 
\bibitem{LedoitWolf20042}  Ledoit, O., Wolf, M. \emph{Honey, I shrunk the sample covariance matrix}.  Journal of Portfolio Management, 30, Volume 4, 2004, 110--119.   
\bibitem{LedoitWolf2012}  Ledoit, O., Wolf, M. \emph{Nonlinear shrinkage estimation of large-dimensional covariance matrices}. Annals of Statistics 40, 2012, 1024--1060.  
\bibitem{RIE0} Ledoit, O.,   P\'ech\'e, S. \emph{Eigenvectors of some large sample covariance matrix ensembles} Probability Theory and Related Fields, 2011, 151.1, 233--264

 \bibitem{Vershynin} Vershynin, R. \emph{High-Dimensional Probability: An Introduction with Applications in Data Science}, Cambridge Series in Statistical and Probabilistic Mathematics, 2018.
 \en{thebibliography}
\en{document}